\documentclass[MSNbibl,number,citesort,dvips]{arxbj}
\usepackage{graphicx}
\aid{0}
\volume{19}
\issue{2}
\pubyear{2013}
\firstpage{676}
\lastpage{719}
\doi{10.3150/11-BEJ412}

\makeatletter
\define@key{bid}{pmid}{}
\define@key{bid}{pmcid}{}

\newtheorem{them}{Theorem}
\newtheorem{prp}{Proposition}
\newtheorem{lem}{Lemma}

\newremark{nota}{Notation}

\newcommand{\argmin}{\arg\min}
\newcommand{\cA}{\mathcal{A}}
\newcommand{\cK}{\mathcal{K}}
\newcommand{\cN}{\mathcal{N}}
\newcommand{\cQ}{\mathcal{Q}}

\newcommand{\bbE}{\mathbb{E}}
\newcommand{\bbH}{\mathbb{H}}
\newcommand{\bbP}{\mathbb{P}}
\newcommand{\bbR}{\mathbb{R}}
\newcommand{\bbV}{\mathbb{V}}
\newcommand{\bbZ}{\mathbb{Z}}

\newcommand{\E}{\mathbb{E}}
\renewcommand{\P}{\mathbb{P}}
\newcommand{\Var}{\operatorname{Var}}

\newcommand{\eps}{\varepsilon}
\renewcommand{\iff}{\Leftrightarrow}
\newcommand{\Stm}{S_m^{(2)}}
\newcommand{\simP}{\sim_{\mathrm P}}
\newcommand\asymP{\asymp_{\mathrm P}}
\renewcommand{\o}{\mathrm{o}}
\renewcommand{\O}{\mathrm{O}}
\newcommand{\OP}{\O_{\mathrm{P}}}

\newcommand{\Vm}{\bbV_m}
\newcommand{\Km}{\cK_m}
\newcommand{\Qm}{\cQ_m}
\newcommand{\Qmt}{\cQ_m^{(t)}}

\newcommand{\ULS}{\operatorname{ULS}}
\newcommand{\ul}{\underline}
\newcommand{\s}{\sigma}

\makeatother

\begin{document}
\begin{frontmatter}

\title{Cluster detection in networks using percolation}
\runtitle{Cluster detection in networks using percolation}

\begin{aug}
\author[a]{\fnms{Ery} \snm{Arias-Castro}\corref{}\thanksref{a}\ead[label=e1]{eariasca@ucsd.edu}} \and
\author[b]{\fnms{Geoffrey R.} \snm{Grimmett}\thanksref{b}\ead[label=e2]{G.R.Grimmett@statslab.cam.ac.uk}}
\runauthor{E. Arias-Castro and G.R. Grimmett} 
\address[a]{Department of Mathematics, University of California, San
Diego, CA 92093-0112, USA.\\ \printead{e1}}
\address[b]{Statistical Laboratory, Centre for Mathematical Sciences,
University of Cambridge, Wilberforce Road, Cambridge CB3 0WB, UK.  \printead{e2}}
\end{aug}

\received{\smonth{6} \syear{2011}}
\revised{\smonth{10} \syear{2011}}

%
\begin{abstract}
We consider the task of detecting a salient cluster in a sensor
network, that is, an undirected graph with a
random variable attached to each node. Motivated by recent research in
environmental statistics and the
drive to compete with the reigning scan statistic, we explore
alternatives based on the percolative properties
of the network. The first method is based on the size of the largest
connected component after removing the
nodes in the network with a value below a given threshold. The second
method is the upper level set
scan test introduced by Patil and Taillie [\textit{Statist. Sci.} \textbf{18}
(2003) 457--465]. We establish the performance
of these methods in an asymptotic decision-
theoretic framework in which the network size increases. These tests
have two advantages over the
more conventional
scan statistic: they do not require previous information about cluster
shape, and they are computationally more
feasible. We make abundant use of percolation theory to derive our
theoretical results, and complement our theory with some numerical experiments.
\end{abstract}
%
%
\end{frontmatter}
%
\section{Introduction}
\label{sec:intro}
We consider the problem of cluster detection in a network. The network
is modeled as a graph, and we
assume that a random variable is observed at each node. The task is
then to detect a cluster, that is, a
connected subset of nodes with values that are larger than usual. There
are a multitude of applications for which this model is relevant;
examples include detection of hazardous
materials (Hills~\cite{hills2001sd}) and target tracking (Li~\textit{et~al.}~\cite{target-sensor})
in sensor networks (Culler, Estrin and Srivastava \cite
{overview-sensor}), and detection of disease
outbreaks (Heffernan \textit{et al.}
\cite{heffernan2004ssp}; Rotz and Hughes~\cite{rotz2004advances};
Wagner \textit{et al.}~\cite{wagner2001esv}).
Pixels in digital images are also sensors, and thus many other
applications are found in the rich literature on
image processing, for example, road tracking (Geman and Jedynak \cite
{road-tracking}) and
fire prevention using satellite imagery (Pozo, Olmo and
Alados-Arboledas~\cite{fire-detection}), and the
detection of cancerous tumors in medical imaging (McInerney and
Terzopoulos~\cite{medical-survey}).

After specifying a distributional model for the observations at the
nodes and a class of clusters to be
detected, the generalized likelihood ratio (GLR) test is the first
method that comes to mind. Indeed, this is by far the
most popular method in practice, and as such, is given different names
in different fields. The likelihood
ratio is known as the scan statistic in spatial statistics (Kulldorff
\cite{Kul,disease-outbreak}) and the corresponding
test as the method of matched filters in engineering (Jain, Zhong and
Dubuisson-Jolly~\cite{deformable-review}; McInerney and Terzopoulos
\cite{medical-survey}). Here we use
the former, where scanning a given cluster $K$ means computing the
likelihood ratio for the simple
alternative where $K$ is the anomalous cluster.
Various forms of scan statistic have been proposed, differing mainly by
the assumptions made on the shape
of the clusters. Most methods assume that the clusters are in some
parametric family (e.g., circular
(Kulldorff and Nagarwalla~\cite{kulldorff1995spatial}), elliptical
(Hobolth, Pedersen and Jensen~\cite{elliptical}; Kulldorff \textit{et al.}
\cite
{kulldorff2006ess})) or, more generally, deformable
templates (Jain, Zhong and Dubuisson-Jolly~\cite{deformable-review}).
Sometimes no explicit shape is
assumed, leading to nonparametric models (Duczmal and
Assun{\c{c}}{\~a}o~\cite{MR2045632}; Kulldorff, Fang and Walsh
\cite{tree-based};
Tango and Takahashi~\cite{tango2005flexibly}).

We consider two alternative nonparametric methods, both based on the
percolative properties of the
network, that is, based on the connected components of the graph after
removing the nodes with values below a given threshold.
The simplest is based on the size of the largest connected component
after thresholding -- the threshold is the only parameter of this method.
If the graph is a one-dimensional lattice, then after thresholding,
this corresponds to the test based on the longest run (Balakrishnan and
Koutras~\cite{MR1882476}), which Chen and Huo~\cite{MR2268049} adapt
for path detection in a thin band.
This test has been studied in a similar context in a series of
papers\footnote{The authors were not aware of this unpublished line of
work until M. Langovoy contacted them in the final stages manuscript
preparation.} (Davies,
Langovoy and Wittich~\cite{langovoy-10}; Langovoy and
 Wittich~\cite{langovoy-11}) under the name of maximum
cluster test.
The idea behind this method is simple. When an anomalous cluster is
indeed present, the values at the nodes belonging to this cluster are
larger than usual and thus more likely to survive the threshold, and
because these nodes are also likely to clump together -- because the
cluster is connected in the graph -- the size of the statistic will be
(stochastically) larger than when no anomalous cluster is present.

More sophisticated, and also parameter-free, is the method based on the
upper level set scan statistic of Patil and Taillie \cite
{MR2109372}, subsequently developed in the context of ecological and
environmental applications (Patil, Joshi
and Koli~\cite{springerlink:10.1007/s10651-010-0140-1}; Patil and
Taillie~\cite{patil-upper}; Patil \textit{et al.}~\cite{patil,MR2297368}).
It is the result of scanning
over the connected components of the graph after thresholding, which is
repeated at all thresholds.
This method obviously is closely related to the scan statistic. It can
be seen as attempting to approximate the scan statistic over all
possible connected components of the graph by restricting the class of
subsets to be scanned to those surviving a threshold. Our results
indicate that this method is in fact more closely related to the
previous one (based on the size of the largest connected component at a
given threshold), and in some sense provides a way to automatically
choose the threshold.

These two percolation-based methods have two significant advantages
over the scan statistic.
First, they do not need to be provided with the shape of the clusters
to be detected. Thus they are valuable
in settings with less previous spatial information. The second
advantage is
computational. The scan statistic tends to be computationally
demanding, even in parametric settings, or
even outright intractable, particularly in nonparametric settings. In
contrast, these two methods are computationally feasible, and their
implementation is fairly straightforward, even for irregular
networks. On the other hand, the scan statistic often relies on the
fast Fourier transform in the square
lattice to scan clusters of known shape over all locations in that network.

In terms of detection performance, we compare these percolation-based
methods to the scan statistic in a
standard asymptotic decision theoretic framework where the network is a
square lattice of growing size and
the variables at the nodes are assumed i.i.d. for nodes inside
(resp., outside) the anomalous cluster.
The performance of the scan statistic in such a framework is well
understood and known to be
(near-) optimal, which makes it the gold standard in detection
(Arias-Castro, Cand\`es and Durand~\cite{cluster}; Arias-Castro,
Donoho and Huo~\cite{MGD}; Perone Pacifico
\textit{et al.}~\cite{perone}; Walther
\cite{MR2604703}).
We find that these two methods are suboptimal for the detection of
hypercubes, an
emblematic parametric class, but are near-optimal for the detection of
self-avoiding paths, an
emblematic nonparametric class.
The main weakness of these percolation-based methods is that when the
per-node signal-to-noise ratio is weak, the connected components after
thresholding are heavily influenced by the whimsical behavior of the
values at the nodes. The scan statistic is very effective in such
situations. Although this rationale seems to apply particularly well in
the case of self-avoiding paths, what makes these methods competitive
in this case is that the problem of detecting such objects is
intrinsically very hard.

The study of the connected components after thresholding is
intrinsically connected to percolation theory (Grimmett \cite
{MR1707339}), an
important branch in probability theory. In fact, when the node values
are i.i.d. -- which is the case when no anomalous cluster is
present -- the only dependence on the distribution at the nodes is the
probability of surviving the threshold, and after thresholding, the
network is a site percolation model. (We introduce and discuss these
notions in detail later in the article.) Our contribution is a careful
analysis of these two nonparametric methods using percolation
theory (Grimmett~\cite{MR1707339}) in a substantial way, thus applying
percolation theory in a sophisticated fashion to shed light on an
important problem in statistics.

The rest of the paper is organized as follows.
In Section~\ref{sec:setting} we formally introduce the framework and
state some
fundamental detection bounds.
In Section~\ref{sec:scan} we describe the standard scan statistic and present
some results on its performance,
showing that it is essentially optimal.
In Section~\ref{sec:cc}, we consider the size of the largest connected
component
after thresholding.
In Section~\ref{sec:uls}, we consider the upper level set scan statistic.
We briefly discuss implementation issues and present some numerical
experiments in Section~\ref{sec:numerics}.
Finally, Section~\ref{sec:discussion} is a discussion section where, in
particular, we mention extensions. We provide proofs in the
\hyperref[appa]{Appendix}.
%
\section{Mathematical framework and fundamental detection bounds}
\label{sec:setting}
For concreteness, and also for its relevance to signal and image processing,
we model the network as a finite subgrid of the regular square lattice
in dimension $d$,
denoted $\Vm:= \{1, \ldots, m\}^d$.
Our analysis is asymptotic in the sense that the network is assumed
to be large, that is, $m \to\infty$. To each node $v \in\Vm$,
we attach a random variable, $X_v$. For example, in the context of
a sensor network, the nodes represent the sensors and the variables
represent the information that they transmit.
The random variables
$\{X_v\dvt v \in\Vm\}$ are assumed to be independent with common distribution
in a certain one-parameter exponential family $\{F_\theta\dvt \theta\in
[0, \theta_\infty)\}$, defined as follows.
Let $\theta_\infty>0$, let $F_0$ be a distribution function with finite
non-zero variance $\sigma_0^2$, and assume the that
moment-generating function $\varphi(\theta):= \int
\mathrm{e}^{x\theta}\,
\mathrm{d}F_0(x)$ is finite for $\theta\in[0, \theta_\infty)$.
Then $F_\theta$ is the distribution function with density $f_\theta(x)
= \exp(\theta x - \log\varphi(\theta))$ with respect to $F_0$.
We assume
further regularity of $F_0$ at later points in this paper.
Note that our results apply to other distributional models as well, as
discussed in Section~\ref{sec:discussion}.

Examples of such a family $\{F_\theta\dvt \theta\in[0,\theta_\infty)\}$
include the following:
\begin{itemize} 
\item\textit{Bernoulli model}: $F_\theta= \operatorname
{Ber}(p_\theta)$,
$p_\theta:= \operatorname{logit}^{-1}(\theta+\theta_0)$, relevant
in sensor
arrays where each sensor transmits one bit (i.e., makes a binary decision)
\item\textit{Poisson model}: $F_\theta= \operatorname{Poi}(\theta
+\theta_0)$,
popular with count data, for example, arising in infectious disease
surveillance systems
\item\textit{Exponential model}: $F_\theta= \operatorname
{Exp}(\theta_0
-\theta
)$ (e.g., to model response times)
\item\textit{Normal location model}: $F_\theta= \cN(\theta+\theta
_0,1)$, standard in signal and image processing, where noise is often
assumed to be Gaussian.
\end{itemize}

\begin{figure}

\includegraphics{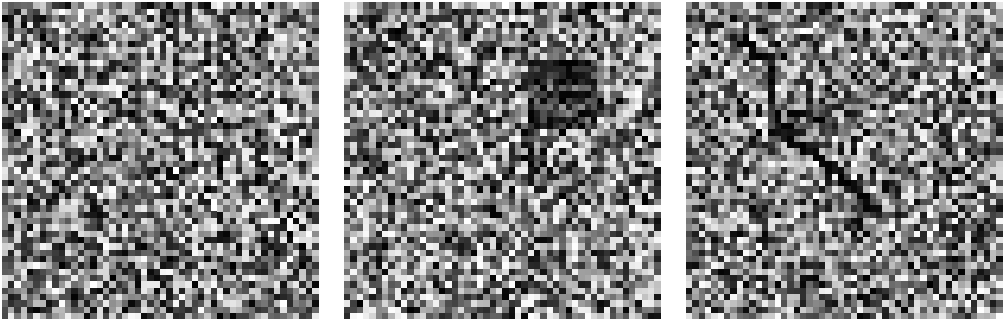}

\caption{This figure illustrates the setting in dimension $d=2$ for a
beta model where $F_0 = \operatorname{Unif}(0,1)$ and $F_\theta=
\operatorname{Beta}(\theta+1,1),
\theta\geq0$. (Left) An instance of the null
hypothesis. (Middle) An instance of an alternative with a square
cluster. (Right) An instance of an alternative with a path.}\label{fig:setting}
\end{figure}

Let $\Km$ be a class of clusters, with a cluster defined as a subset of
nodes connected in the graph.
Under the null hypothesis, all of the variables at the nodes have
distribution $F_{0}$, that is,
\[
\bbH^m_0 \dvt X_v \sim F_{0}\qquad \forall v \in\Vm.
\]
Under the particular alternative where $K \in\Km$ is anomalous, the
variables indexed by $K$ have distribution $F_{\theta_m}$ for some
$\theta_m > 0$, that is,
\[
\bbH^m_{1,K}\dvt X_v \sim F_{\theta_m}\qquad \forall v \in K;\qquad X_v
\sim F_{0}\qquad \forall v \notin K.
\]
We are interested in the situation where the anomalous cluster $K$ is
unknown, namely in testing $\bbH^m_0$ against $\bbH^m_1 := \bigcup_{K
\in\Km} \bbH^m_{1,K}$. We illustrate the setting in Figure~\ref{fig:setting}
in the context of the two-dimensional square grid.

Let $\Km$ denote a cluster class for $\Vm$.
As usual, a test $T$ is a function of the data, $T = T(X_v\dvt v \in\Vm
)$, that takes values in $\{0,1\}$, with $T = 1$ corresponding to a
rejection, meaning a decision in favor of $\bbH^m_1$.
For a test $T$, we define its worst-case risk as the sum of its
probability of type I error and its probability of type II maximized
over the anomalous clusters in the class
\[
\gamma_m(T) = \P(T = 1|\bbH^m_0) + \max_{K \in\Km} \P(T = 0|\bbH
^m_{1,K}).
\]
A method is formally defined as a sequence of tests $(T_m)$ for testing
$\bbH^m_0$ versus $\bbH^m_1$.
We say that a method $(T_m)$ is asymptotically powerless if
\[
\liminf_{m \to\infty} \gamma_m(T_m) \geq1.
\]
This amounts to saying that as the size of the network increases, the
method $(T_m)$ is not substantially better than random guessing.
Conversely, a method $(T_m)$ is asymptotically powerful if
\[
\lim_{m \to\infty} \gamma_{m}(T_m) = 0.
\]
%
The minimax risk is defined as $\gamma_m^* := \inf_T \gamma_m(T)$, and
we say that a method is $(T_m)$ (asymptotically) optimal if $\gamma
_{m}(T_m) \to0$ whenever $\gamma_m^* \to0$. Everything else fixed,
the latter depends on the behavior of $\theta_m$ when $m$ becomes
large. We say that $(T_m)$ is optimal up to a multiplicative constant
$C \geq1$ if $\gamma_{m}(T_m) \to0$ under $C \theta_m$ whenever
$\gamma_m^* \to0$ under $\theta_m$. We say that $(T_m)$ is
near-optimal if the same is true with $C$ replaced by $C_m \to\infty$
with $\log C_m = \mathrm{o}(\log\theta_m)$. (This occurs here only when
$\theta
_m \to0$ polynomially fast and $C_m \to\infty$ poly-logarithmically fast.)

We focus on situations where the clusters in the class $\Km$ are of same
size, increasing with $m$ but negligible compared with the size of the
entire network. We do so for the sake of simplicity; more general
results could be obtained as in Arias-Castro, Cand\`es and Durand \cite
{cluster}, Arias-Castro,
Donoho and Huo~\cite{MGD}, Perone Pacifico
\textit{et al.}~\cite{perone}, Walther~\cite{MR2604703}
without additional difficulty. Assuming a large
anomalous cluster allows us to state general results
applying to a wide range of one-parameter exponential families
(via the central limit theorem). In addition, note that
on the one hand, reliably detecting a cluster of bounded size is
impossible in the Bernoulli model or any other model where $F_0$
has finite support, whereas on the other hand, detecting a cluster of size
comparable to that of the entire network is in some sense trivial,
given that
the simple test based on the total sum $\sum_{v \in\Vm} X_v$ is
optimal up to a multiplicative constant.
%

We consider two emblematic classes of clusters, in some sense at the
opposite extremes:
\begin{itemize}
\item
\textit{Hypercube detection.} Let $\Km$ denote the class of \textit
{hypercubes} within $\Vm$ of sidelength $[m^\alpha]$ with $0 < \alpha<
1$. This class is parametric, with the location of the hypercube the
only parameter.
\item
\textit{Path detection.} Let $\Km$ denote the class of \textit{loopless
paths} within $\Vm$ of length $[m^\alpha]$ with $0 < \alpha< 1$. This
class is nonparametric, in the sense that its cardinality is
exponential in the length of the paths.
\end{itemize}
See Figure~\ref{fig:setting} for an illustration. (Note that a
hypercube of
side length $k$ may be seen as a loopless path of length $k^d$.)
Although we obtain results for both, our main focus is in the setting
of hypercube detection, which is relevant to a wider range of
applications, in fact any situation where the task is to detect a shape
that is not filamentary. The situation exemplified in the setting of
path detection may be relevant in target tracking from video, or the
detection of cracks
in materials in non-destructive testing.
Note that the two settings coincide in dimension one.

We state fundamental detection bounds for each setting.
The following result is standard (see, e.g., Arias-Castro, Cand\`es and
Durand~\cite{cluster}; Arias-Castro, Donoho and Huo~\cite{MGD}).
Remember that $\sigma_0^2$ denotes the variance of $F_0$.
\begin{lem} \label{lem:cube}
In hypercube detection, all methods are asymptotically powerless if
\[
\limsup_{m \to\infty} (\log m)^{-1/2} m^{d\alpha/2} \theta_m <
\sigma_0 \sqrt{2 d (1 -\alpha)}.
\]
\end{lem}

In fact, the conclusions of Lemma~\ref{lem:cube} apply for a wide
variety of
parametric classes, such as discs, a popular model in disease outbreak
detection (Kulldorff and Nagarwalla~\cite{kulldorff1995spatial}), as
well as to nonparametric
classes of blob-like clusters (see Arias-Castro, Cand\`es and Durand
\cite{cluster}; Arias-Castro, Donoho and Huo~\cite{MGD}).

The following result is taken from Arias-Castro \textit{et al.}~\cite{MR2435454}.
\begin{lem} \label{lem:path}
In path detection, all methods are asymptotically powerless if $\lim_{m
\to\infty} \theta_m \*(\log m)(\log\log m)^{1/2} = 0$, in dimension
$d=2$, and the same is true in dimension $d \geq3$ if $\limsup_{m \to
\infty} \theta_m < \theta_*$, where $\theta_* > 0$ depends only on $d$.
\end{lem}

In dimension $d \geq4$, $\theta_*$ may be taken to be the unique
solution to
\[
\rho\varphi(2\theta) - \varphi(\theta)^2 = 0,
\]
where $\rho$ is the return probability of a symmetric random walk in
dimension $d$.
%
%
\section{The scan statistic}
\label{sec:scan}
For a subset of nodes $K \subset\bbV$, let $|K|$ denote its size and define
\[
\bar{X}_K = \frac{1}{|K|} \sum_{v \in K} X_v.
\]
Given a cluster class $\cK$, we define the (simple) scan statistic as
%
\begin{equation}\label{scan}
\max_{K \in\cK} \sqrt{|K|} (\bar{X}_K -\mu_0),
\end{equation}
where $\mu_0$ is the mean of $F_0$.
If $\mu_0$ is not available, we may use the grand mean $\bar{X}_{\Vm
}$ instead.
In Appendix \hyperref[sec:proof-scan]{B}, we derive this form of the scan
statistic as an approximation to the scan statistic of Kulldorff
\cite{Kul}, which is, strictly speaking, the GLR and arguably the most
popular version, particularly in spatial statistics.
We use this simpler form to streamline our theoretical analysis.

The test that rejects for large values of the scan statistic (\ref
{scan}), which we call the scan test, is near-optimal in a wide range of
settings (Arias-Castro, Cand\`es and Durand~\cite{cluster}; Arias-Castro, Donoho and Huo \cite
{MGD}; Walther~\cite{MR2604703}). In
particular, in the context
of a class of hypercubes, and in fact many other parametric classes,
this test is asymptotically optimal to the exact multiplicative constant.
\begin{lem} \label{lem:scan-cube}
In hypercube detection, the scan test is asymptotically powerful if
\[
\liminf_{m \to\infty} (\log m)^{-1/2} m^{d\alpha/2} \theta_m >
\sigma_0 \sqrt{2 d (1 -\alpha)}.
\]
\end{lem}

In the context of a class of paths, the following result states that
the scan test detects if $\theta_m$ is bounded away from 0 and
sufficiently large. Note that this does not match the order of
magnitude of the lower bound given in dimension $d=2$.
Let $\Lambda(\theta) = \log\varphi(\theta)$ and
$
\Lambda^*(x) = \sup_{\theta\geq0} [\theta x - \Lambda(\theta
) ].
$
($\Lambda^*$ is the rate function of $F_0$ when $x \geq\mu_0$.)
The following result is established in Arias-Castro \textit{et al.}~\cite{MR2435454}.
\begin{lem} \label{lem:scan-path}
In path detection, the scan test is asymptotically powerful if
\[
\liminf_{m \to\infty} \theta_m > \theta_* := (\Lambda^* \circ
\Lambda
')^{-1}(\log(2d)).
\]
\end{lem}
%
\section{Size of the largest open cluster}
\label{sec:cc}
We study the test based on the size of the largest connected component
after thresholding the values at the nodes. This test was independently
considered in a series of papers (Davies, Langovoy and Wittich~\cite{langovoy-10}; Langovoy and Wittich \cite
{langovoy-11}). Our
results are seen to sharpen and elaborate on these results. In
particular, we study this test under all three regimes (subcritical,
supercritical, and critical).

Adapting terminology from percolation theory (Grimmett \cite
{MR1707339}), for a
threshold $t \in\bbR$, we say that a subset $K \subset\bbV$ is open
(at threshold $t$) if $X_v > t$ for all $v \in K$. Let $S_m(t)$
(resp., $S_K(t)$) denote the size of the largest open cluster within
$\Vm
$ (resp., within $K$). The analysis of the test based on $S_m(t)$, which
we call the largest open cluster (LOC) test, boils down to bounding the
size of $S_m(t)$ from above, under $\bbH^m_0$, and, because $S_m(t)
\geq S_K(t)$, bounding the size of $S_K(t)$ from below, under $\bbH
^m_{1,K}$. Define $\xi_v(t) = \mathbf{I}\{X_v > t\}$, which is Bernoulli
with parameter $p_\theta(t) := \bbP_\theta(X_v > t)$. The process
$(\xi
_v(t)\dvt v \in\bbV_m)$ is a site percolation model (Grimmett~\cite{MR1707339}).
In general, consider a process $(\xi_v\dvt v \in\Vm)$ i.i.d. Bernoulli
with parameter $p$, and let $S_m$ denote the size of the largest open
cluster within $\Vm$. In dimension $d = 1$, this process may be seen as
a sequence of coin tosses, and $S_m$ viewed as the longest heads run in
that sequence. In this context, the Erd\H os--R\'enyi Law (Erd{\H{o}}s
and R{\'e}nyi~\cite{MR0272026}) says that
%
\begin{equation}\label{S1-conv}
\frac{S_m}{\log m} \to\frac{1}{\log(1/p)}, \qquad\mbox{almost surely}.
\end{equation}
In higher dimensions $d \geq2$, the situation is much more involved.
Let $p_c$ denote the critical probability for site percolation in $\bbZ
^d$, defined as the supremum over all $p \in(0,1)$ such that the size
of the open cluster at the origin, denoted by $S$, is finite with
probability 1. (The dependency in $d$ is left implicit.) We consider
the subcritical ($p_0(t) < p_c$), supercritical ($p_0(t) > p_c$), and
near-critical ($p_0(t) \approx p_c$) cases separately.
%
\subsection{Subcritical percolation}
\label{sec:cc-sub}
In the subcritical case, where $t$ is such that $p_0(t) < p_c$, we are
able to obtain precise, rigorous results on the performance of the test
based on $S_m(t)$ in terms of the function $\zeta_p$, implicitly
defined as
%
\begin{equation}\label{zeta}
\zeta_p := - \lim_{k \to\infty} \frac{1}{k} \log\mathbb{P}(S
\geq k) = - \lim
_{k \to\infty} \frac{1}{k} \log\mathbb{P}(S = k)
\end{equation}
(see Grimmett~\cite{MR1707339}, Section 6.3). Again, the dependency in $d$
is left
implicit.
As a function of $p \in(0,p_c)$, $\zeta_p$ is continuous and strictly
decreasing, with limits $\infty$ at $p=0$ and 0 at $p=p_c$ (see Lemma
\ref{lem:zeta}), whereas $\zeta_p = 0$ for $p \geq p_c$.
In the \hyperref[appa]{Appendix}, we include a proof that
%
\begin{equation}\label{S-conv}
\frac{S_m}{\log m} \to\frac{d}{\zeta_p}, \qquad\mbox{in probability}
\end{equation}
for a subcritical threshold $p < p_c$.

The convergence result in (\ref{S-conv}) may be used to bound $S_m(t)$
under the null by taking $p = p_0(t)$.
Under the alternative, if we consider a class of hypercubes, then (\ref
{S-conv}) also may be used to bound $S_K(t)$, because $K$ is a scaled
version of $\Vm$.
\begin{them} \label{thm:cc-cube}
In hypercube detection, the test based on $S_m(t)$, with $t$ fixed such
that $0 < p_0(t) < p_c$, is asymptotically powerful if
$\liminf_{m \to\infty} \theta_m > \theta_*(t)$,
and asymptotically powerless if
$\limsup_{m \to\infty} \theta_m < \theta_*(t)$,
where $\theta_*(t)$ is the unique solution to $\zeta_{p_\theta(t)} =
\alpha\zeta_{p_0(t)}$.
\end{them}

Note that when $t$ is fixed, $\zeta_{p_\theta(t)}$ as a function of
$\theta$ is continuous and strictly
strictly decreasing, by the fact that $p_\theta(t)$ is continuous and strictly
increasing in $\theta$ (Brown~\cite{MR882001}, Cor. 2.6, 2.22) and
$\zeta_p$ is
continuous
and strictly decreasing in $p$ (Lemma~\ref{lem:zeta}).
Therefore, $\theta_*(t)$ in the theorem is well defined.

If instead, we consider a class of paths, then (\ref{S1-conv}) may be
used to bound $S_K(t)$, because $K$ is a scaled version of the lattice
in dimension 1. In congruence with (\ref{S1-conv}), we define $\zeta
^1_p = \log(1/p)$.
\begin{them} \label{thm:cc-path}
In path detection, the test based on $S_m(t)$, with $t$ fixed such that
$0 < p_0(t) < p_c$, is asymptotically powerful if
$\liminf_{m \to\infty} \theta_m > \theta_*^+(t)$,
and asymptotically powerless if
$\limsup_{m \to\infty} \theta_m < \theta_*^-(t)$,
where $\theta_*^+(t)$ (resp., $\theta_*^-(t)$) is the unique solution to
$d \zeta^1_{p_\theta(t)} = \alpha\zeta_{p_0(t)}$ (resp., $d \zeta
^1_{p_\theta(t)} = \alpha\zeta_{p_0(t)}$).
\end{them}

Note that in dimension $d \geq2$, the result is not sharp, because we
always have $\theta_*^+(t) > \theta_*^-(t)$.
We believe that sharper forms of this result may be substantially more
involved, and for this reason we have not pursued this.

Qualitatively, the message is that for both hypercube detection and
path detection, the subcritical LOC test requires that $\theta_m$ be
larger than a constant to be effective. Compared with the scan
statistic, this makes it grossly suboptimal when detecting hypercubes
and comparable (up to a multiplicative constant in $\theta_m$) when
detecting self-avoiding paths.

What if we let $t = t_m \to\infty$, so that $p_0(t_m) \to0$? Then the
test based on $S_m(t_m)$ is powerless under some additional conditions
on $F_0$.
For $b, C \geq0$, consider the following class of approximately
exponential power ($\operatorname{AEP}$) distributions, sometimes
called Subbotin distributions:
\[
\operatorname{AEP}(b, C) = \{F\dvt x^{-b} \log\bar{F}(x) \to-C, x \to
\infty\}.
\]
($\bar{F}(x) := 1 -F(x)$ is the survival distribution function of $X
\sim F$.)
For example, $\operatorname{Exp}(\lambda) \in\operatorname{AEP}(1,
\lambda)$ and $\cN(\mu
, \sigma^2) \in\operatorname{AEP}(2, 1/(2\sigma^2))$, whereas
$\operatorname{Poi}(\lambda
)$ behaves roughly as a distribution in $\operatorname{AEP}(1, C)$.
%
\begin{prp} \label{prp:t-large}
Assume that $F_0 \in\operatorname{AEP}(b, C)$ for some $b > 1$ and $C >
0$. In
hypercube detection, the test based on $S_m(t)$ is asymptotically
powerless when $t = t_m \to\infty$, unless $\theta_m \to\infty$.
\end{prp}
%
\subsection{Supercritical percolation}
\label{sec:cc-cube-sup}
Here we consider the supercritical regime, where $p_0(t) > p_c$. (Note
that necessarily $d \geq2$ for $p_c = 1$ in dimension 1.) In this
setting, too, the size of the largest cluster is well understood.
Let $\Theta_p$ be the probability that the open cluster at the origin
is infinite, and note that $\Theta_p>0$
for $p>p_c$, by the definition of $p_c$. We have with probability 1 that
\[
\frac{S_m}{|\Vm|} \to\Theta_p
\]
(see Falconer and Grimmett~\cite{FG92}, Lemma 2 and proof, Penrose and
Pisztora~\cite{MR1372330}, Theorem~4,
Pisztora~\cite{Pisz96}). In fact
(with probability $1-\o(1)$), the largest open cluster within $\Vm$
is unique,
and the foregoing statement says that it occupies a non-negligible
fraction of $\Vm$.
With a supercritical choice of threshold, the LOC test is powerless for
any $\theta$ if the anomalous cluster is too small, specifically if
$\alpha< 1/2$ in the setting of hypercube detection. Indeed, we have
the following result.
\begin{them} \label{thm:cc-cube-sup}
In hypercube detection, the test based on $S_m(t)$, with $t$ fixed such
that $p_c < p_0(t) < 1$, is asymptotically powerful if $\alpha\geq
1/2$ and
$\lim_{m \to\infty} \theta_m m^{(\alpha-1/2) d} = \infty$,
and asymptotically powerless if $\alpha< 1/2$ or if
$\lim_{m \to\infty} \theta_m m^{(\alpha-1/2) d} = 0$.
\end{them}
%

Thus, for the detection of small clusters, a supercritical LOC test is
potentially worthless, whereas for larger clusters it improves
substantially on the performance of a subcritical LOC test, although it
is still suboptimal compared with the scan statistic. (Indeed,
comparing the exponents when $\alpha\geq1/2$, we have $(\alpha-1/2)
d < \alpha d/2$, because $\alpha< 1$.)
We mention that in the context of path detection, the same arguments
show that the LOC test for any choice of supercritical threshold is
asymptotically powerless.
%
\subsection{Critical percolation}
\label{sec:cc-cri}
If our goal is to choose a threshold $t$ so as to maximize the
difference in size for the largest open cluster under the null and
under an alternative, then we are necessarily in the neighborhood
of the percolation phase transition, which is to say that $|p-p_c|$ is
small. (Again, here we assume $d \geq2$.) The percolation model is
not fully understood in the critical regime, which poses a serious
obstacle to a rigorous statistical analysis.
(See Grimmett~\cite{MR1707339}, Chapter 9, for a general discussion of this
percolation regime.)
We base our discussion on the work of Borgs \textit{et al.} \cite
{MR1868996}. Let $\pi_m(p)$ denote the probability that the open
cluster at the origin reaches outside the box $[-m,m]^d$, and let $\xi
(p)$ denote the correlation length, defined as
\[
\frac{1}{\xi(p)} := -\lim_{m \to\infty} \frac{1}{m} \log\pi_m(p).
\]
%
Note that, with $\xi$ thus defined, $\xi(p)<\infty$ if and only if
$p <
p_c$. The critical exponent for (subcritical) correlation
length is postulated to be
\[
\nu:= -\lim_{p \nearrow p_c} \frac{\log\xi(p)}{\log|p - p_c|}.
\]
It is not known whether the limit exists for all dimensions, but
it is known that $0 < \nu< \infty$ whenever it exists.
It is shown in Borgs \textit{et al.}~\cite{MR1868996} that, subject to the
existence of this
limit together with other
scaling assumptions, when $p = p_m$ varies with $m$,
%
\begin{equation}\label{S-cri}
S_m \asymP
\cases{\log m,& \quad$\mbox{if, for some }\nu'>\nu,  m^{1/\nu'} (p_m
-p_c) \to
-\infty$,\cr
m^d, &\quad$\mbox{if, for some } \nu'>\nu,  m^{1/\nu'} (p_m -p_c) \to
\infty$,}
\end{equation}
where $X_m \asymP Y_m$ means that there exists a constant $C \in(0,
\infty)$
such that $C^{-1} \leq X_m/Y_m \leq C$ in probability.
The scaling assumptions of Borgs \textit{et al.}~\cite{MR1868996} are believed
to hold if
and only if the number
$d$ of dimensions satisfies $2\le d \le6$, and they are proved for
$d=2$. The work of Borgs \textit{et al.}~\cite{MR1868996}
was directed at bond percolation only, but similar results are expected
for site percolation.

It is known that $\nu=4/3$ for site percolation on the triangular
lattice (see Smirnov and Werner~\cite{SmWe}), and it is believed that
this holds for
percolation on any two-dimensional lattice. As described in Grimmett
\cite{MR1707339}, Section 10.4,
it is believed that $\nu=1/2$ for $d \ge6$, and this has been proved
for $d \ge19$ and for the so-called ``spread-out model'' in $7$ and more
dimensions (Hara, van der Hofstad and Slade~\cite{MR1959796}).

Subject to the assumption that (\ref{S-cri}) holds,
we establish the power of the test based on $S_m(t)$ when choosing $t =
t_m$ near criticality.
We assume that there exists $t_c$ such that $p_0(t_c)=p_c$, and that
$p_0(t)$ is a continuous function of $t$ in a neighborhood of $t_c$.
%
%
%
\begin{them} \label{thm:cc-cube-cri}
Let $t_m \geq t_c$ be such that $p_c -p_0(t_m) \asymp m^{-1/\nu'}$ for
some $\nu'>\nu$.
In hypercube detection, assuming that (\ref{S-cri}) holds, the test
based on $S_m(t_m)$ is asymptotically powerful if $\liminf_{m\to
\infty
} \theta_m m^{\alpha/\nu'}$ is sufficiently large.
\end{them}
%
%

Compared with a subcritical choice of threshold, which requires that
$\theta_m$ be bounded away from 0 for the test to have any power, as
seen in Theorem~\ref{thm:cc-cube}, with a near-critical choice of
threshold, the
test is able to detect at polynomially small $\theta_m$. In particular,
with a proper choice of threshold, the test is powerful for $\theta_m$
of order $m^{-\alpha/\nu'}$ with $\nu' > \nu$. Note that, by Lemma
\ref{lem:cube}, all methods are asymptotically powerless if $\theta
_m$ is of
order $m^{-d\alpha/2}$, implying that $\alpha/\nu\le d\alpha/2$. We
thus obtain the inequality $\nu\ge2/d$.
This may be compared with the scaling relation (Grimmett \cite
{MR1707339}, Equation~(9.23)) stating that $d\nu=2-a$, where
$a$ ($<0$) is the percolation critical exponent for the number of
clusters per vertex. It is believed that
$a=-\frac23$ when $d=2$ and $a=-1$ when $d \ge6$.
Compared with the performance at supercriticality, the test at
near-criticality (with a proper choice of threshold) is superior if
$(\alpha-\frac12) d < \alpha/\nu$, which is equivalent to $\alpha<
(1-a/2)/(1-a)$. For example, with $a=-\frac23$, the near-critical LOC
test is superior when $\alpha< \frac34$.
%
%
\section{The upper level set scan statistic}
\label{sec:uls}
For a threshold $t$, let $\Qmt$ denote the (random) class of clusters
within $\Vm$ open at $t$, and let $\Qm^* = \bigcup_{t} \Qmt$, which is
also random. Patil and Taillie~\cite{MR2109372} suggested scanning
the clusters in $\Qm^*$. To facilitate a rigorous mathematical analysis
of its performance, we consider the upper level set ($\operatorname{ULS}$) scan at a
given threshold $t$, and use the simple scan described in Section \ref
{sec:scan}. Specifically, in correspondence with (\ref{scan}), we define
the (simple) $\operatorname{ULS}$ scan statistic at threshold $t$ as
%
\begin{equation}\label{uls}
U_m(t, k_m) = \max\bigl\{ \sqrt{|K|} (\bar{X}_K -\mu_{0|t}) \dvt K \in
\Qmt, |K| \geq k_m \bigr\},
\end{equation}
where $\mu_{0|t}$ (resp., $\sigma^2_{0|t}$) is the the mean
(resp., variance) of $X_v|X_v > t$ when $X_v \sim F_0$,
and $(k_m)$ is a non-decreasing sequence of positive integers.
The $\operatorname{ULS}$ scan statistic of Patil and Taillie~\cite{MR2109372}
corresponds (in its simple form) to
%
\begin{equation}\label{uls-orig}
\ULS_m = \max_{t \in\bbR} \frac{U_m(t, 1)}{\sigma_{0|t}}.
\end{equation}
If $\mu_{0|t}$ and/or $\sigma_{0|t}^2$ are not available, we may use
their empirical versions based on the $X_v$ that survive the threshold
$t$. We restrict the scan to clusters of size at least $k_m$
to increase power, because the behavior of $U_m(t)$ is, as we show
later, completely driven by the smallest open clusters that are
scanned, at least when $t$ is subcritical.
%
We present the rest of our discussion in terms of subcritical,
supercritical, and near-critical choices of threshold. We then conclude
with a result on the performance of the $\operatorname{ULS}$ scan test across all thresholds.
%
\subsection{Subcritical threshold}
\label{sec:uls-sub}
We start by describing the behavior of $U_m(t, k_m)$ under the null.
Let $F_{\theta|t}$ denote the distribution of $X_v|X_v > t$ under
$F_\theta$, and let $\mu_{\theta|t}$ and $\Lambda^*_{\theta|t}$ denote
its mean and rate function, respectively.
Also, when $0 < \beta< 1/\zeta_{p_\theta(t)}$, or $\beta= 0$ and $F_0
\in\operatorname{AEP}(b, C)$ for some $b \geq2$ and $C > 0$, let
$\gamma
_{\theta|t}(\beta) := \gamma(F_{\theta|t}, \mu_{0|t}, \zeta
_{p_\theta
(t)}, \beta)$, where $\gamma$ is the function defined in Lemma \ref
{lem:gamma}.
Note that $\gamma_{\theta|t}(\beta)$ can be computed explicitly in some
cases, like the normal location model, and $\gamma_{\theta|t}(\beta)
\sim(\mu_{\theta|t} -\mu_{0|t})^2/\zeta_{p_\theta(t)}$ when
$\theta
\nearrow\theta_c(t)$, defined (when it exists) as the solution to
$p_{\theta}(t) = p_c$.
%
%
%
%
\begin{lem} \label{lem:uls-sub}
Assume that $\theta\geq0$ and $t$ is fixed such that $0 < p_\theta(t)
< p_c$ and that $k_m/\log m \to d \beta$ for some $\beta\geq0$. Then,
under $F_\theta$ on $\Vm$, the following holds in probability:
\begin{enumerate}
\item If $\beta> 1/\zeta_{p_\theta(t)}$, then $U_m(t, k_m) = 0$ for
$m$ large enough.
\item If $0 < \beta< 1/\zeta_{p_\theta(t)}$, then
\[
(\log m)^{-1/2} U_m(t, k_m) \to(d \gamma_{\theta|t}(\beta))^{1/2}.
\]
\item If $\beta= 0$ and $F_0 \in\operatorname{AEP}(b, C)$ for some $b
\geq1$
and $C > 0$, then
\begin{enumerate}[(b)]
\item[(a)] If $b \geq2$, the convergence in Part 2 applies;
\item[(b)]If $b < 2$,
\[
k_m^{1/b -1/2} (\log m)^{-1/b} U_m(t, k_m) \to(d/C)^{1/b}.
\]
%
\end{enumerate}
\end{enumerate}
\end{lem}

In the last case, where $\beta= 0$, the behavior of $U_m(t)$ is
influenced by the very large deviations of $F_{\theta|t}^{*k}$ for $k
\geq k_m$. (The symbol $*$ denotes convolution.) We choose to state a
result for $\operatorname{AEP}$ distributions, for which the very large
deviations
resemble the large deviations.

Based on Lemma~\ref{lem:uls-sub}, we establish the performance of the
$\operatorname{ULS}$ scan
statistic. We start by arguing that choosing $k_m$ such that $k_m/\log
m \to0$ leads to a test that may potentially have less power than the
test based on the largest cluster after thresholding. Indeed, the
behavior of the $\operatorname{ULS}$ scan statistic does not depend on $\theta$ as long
as $\theta< \theta_c(t)$.
\begin{prp} \label{prp:uls-bad}
Assume that $F_0 \in\operatorname{AEP}(b, C)$ for some $b \in(1,2)$ and
$C >
0$. In hypercube detection, the test based on $U_m(t, k_m)$, with
$t$ fixed such that $0 < p_0(t) < p_c$ and $k_m/\log m \to0$, is
asymptotically powerless if $\limsup_{m \to\infty} \theta_m <
\theta_c(t)$.
\end{prp}

For example, in the setting just described with $d = 1$, the $\operatorname{ULS}$ scan
test has (asymptotically) no power unless $\theta_m \to\infty$,
whereas the test based on the size of the largest cluster after
thresholding is, by Theorem~\ref{thm:cc-cube}, asymptotically powerful if
$\liminf_{m \to\infty} \theta_m$ is large enough.
We therefore choose a sequence $k_m$ comparable in magnitude to $\log
m$ and state the performance of the $\operatorname{ULS}$ scan test in this case.
\begin{them} \label{thm:uls-cube}
In hypercube detection, the test based on $U_m(t, k_m)$, with $t$
fixed such that $0 < p_0(t) < p_c$ and $k_m/\log m \to d \beta$ with $0
< \beta< 1/\zeta_{p_0(t)}$, is asymptotically powerful if
$\liminf_{m \to\infty} \theta_m > \theta_*(t)$
and asymptotically powerless if
$\limsup_{m \to\infty} \theta_m < \theta_*(t)$,
where $\theta_*(t)$ is the unique solution to $\alpha\gamma_{\theta
|t}(\beta) = \gamma_{0|t}(\beta)$.
\end{them}

Note that $\theta_*(t)$ is well defined by Lemma~\ref{lem:gamma-prop}
and that
$\theta_*(t) < \theta_c$ as long as $\alpha> 0$.
In any case, the test based on $U_m(t, k_m)$ with a subcritical
threshold $t$ is, in the setting of hypercube detection, asymptotically
powerless when $\theta_m \to0$, just like the LOC test. In essence,
the two tests are qualitatively comparable in this setting. This is
also true in the context of path detection. Let $\gamma^1_{\theta
|t}(\beta)$ denote $\gamma_{\theta|t}(\beta)$ in dimension 1.
\begin{them} \label{thm:uls-path}
In path detection, the test based on $U_m(t, k_m)$, with $t$ fixed
such that $0 < p_0(t) < p_c$ and $k_m/\log m \to d \beta$ with $0 <
\beta< 1/\zeta_{p_0(t)}$, is asymptotically powerful if
$\liminf_{m \to\infty} \theta_m > \theta_*^+(t)$,
and asymptotically powerless if
$\limsup_{m \to\infty} \theta_m < \theta_*^-(t)$,
where $\theta_*^+(t)$ (resp., $\theta_*^-(t)$) is the unique solution to
$\alpha\gamma^1_{\theta|t}(\beta) = \gamma_{0|t}(\beta)$
(resp., $\alpha\gamma_{\theta|t}(\beta) = \gamma_{0|t}(\beta)$).
\end{them}

As in Theorem~\ref{thm:cc-path}, the result is not as sharp.

Qualitatively, we see that the performance of the subcritical $\operatorname{ULS}$ scan
and LOC tests are comparable for both hypercube detection and path detection.
%
\subsection{Supercritical threshold}
\label{sec:uls-sup}
Here we consider the choice of a supercritical threshold, where $t$ is
fixed such that $p_0(t) > p_c$. We already saw in Section~\ref{sec:cc-cube-sup}
that the largest open cluster is unique and occupies a non-negligible
fraction of the entire network. This is actually true both under the
null and under an alternative.
The $\operatorname{ULS}$ scan test based solely on the largest open cluster is
comparable to the test based on the grand mean after thresholding. In
turn, assuming $t$ is fixed, this test is asymptotically powerful when
$m^{(\alpha-1/2) d} \theta_m \to\infty$, and asymptotically powerless
if $\alpha\leq1/2$ and $\theta_m$ is bounded. (This is easily seen
using Chebyshev's inequality.) This is comparable to the LOC test at
supercriticality.

In general, the $\operatorname{ULS}$ scan statistic includes other (smaller) open clusters.
The story of the second-largest cluster of supercritical percolation in
a box is not yet complete, and for this reason the behavior of the $\operatorname{ULS}$
scan statistic remains incompletely understood. The difficulty arises
from the possibility that the second-largest cluster in $\Vm$ might lie
at its boundary. Whether or not this occurs depends on the outcome of a
calculation (yet to be done) of energy/entropy type involving so-called
``droplets'' near the boundary of $\Vm$ (see, e.g., Bodineau, Ioffe and
Velenik~\cite{winterb}).
To simplify the discussion, we finesse this problem by working where
necessary on $\Vm$ with \emph{toroidal boundary conditions}. That is,
whenever we make statements concerning supercritical percolation on the
graph $\Vm$, we may add edges connecting sites on its boundary as
follows: when $d=2$, for $k=1,2,\ldots,m$, an additional edge is placed
between site $(1,k)$ and site $(m,k)$, and similarly between $(k,1)$
and $(k,m)$.

In proving exact asymptotics for test statistics under the null, we
assume toroidal boundary conditions. Our results on asymptotic power do
not require such exact results but require only orders of magnitude,
which do not need the toroidal assumption. We emphasize that similar
results are expected to hold with ``free'' (i.e., without the extra
edges) rather than
toroidal boundary conditions. Once the percolation picture is better
understood, such results
will follow in the same manner as those presented in this paper. Our results
for the torus are also valid if instead we discount open clusters that
touch the boundary of $\Vm$. Details of this are omitted, and the
proofs are essentially the same.

When working on the torus, the second-largest cluster is controlled
through the following calculation. Cerf~\cite{Cerf06} proved that the limit
%
\begin{equation}\label{delta}
\delta_p := - \lim_{k \to\infty} k^{-(d-1)/d} \log\mathbb
{P}(\infty> S \geq k) = - \lim_{k \to\infty} k^{-(d-1)/d} \log
\mathbb{P}(S = k),
\end{equation}
exists, with $0 < \delta_p < \infty$ for all fixed $p \in(p_c,1)$. The
dependency on $d$ is left implicit.

A result similar to Lemma~\ref{lem:uls-sub} holds with $\delta_p$
playing the
role of $\zeta_p$ and the exponent of $\log m$ changed in places. It
turns out that we need this result only when $\theta= 0$.
For $\beta> 0$ and a supercritical $t$, let $\gamma_{0|t}(\beta) :=
\gamma(F_{0|t}, \mu_{0|t}, 0, \beta)$, defined in Lemma~\ref{lem:gamma}.
\begin{lem} \label{lem:uls-sup}
Assume that $t$ is fixed such that $p_c < p_0(t) < 1$ and that
$k_m/\log m \to d \beta$ and $k_m^{(d-1)/d}/\log m \to d \beta'$ for
some $0 \leq\beta, \beta' \leq\infty$. Then, under the null, the
following holds in probability
on the torus $\Vm$:
\begin{enumerate}
\item If $\beta' > 1/\delta_{p_0(t)}$, then $U_m(t, k_m) = \O(1)$.
\item If $0 \leq\beta' < 1/\delta_{p_0(t)}$ and $\beta= \infty$, then
\[
(\log m)^{-1/2} U_m(t, k_m) \to\sigma_{0|t} \bigl[2 d \bigl(1 -\beta' \delta
_{p_0(t)}\bigr)\bigr]^{1/2},
\]
where $\sigma_{0|t}^2 := \operatorname{Var}(F_{0|t})$.
\item If $\beta< \infty$, then the conclusions of Lemma~\ref{lem:uls-sub}
apply. (Note that $\zeta_{p_0(t)} = 0$.)
\end{enumerate}
\end{lem}

Based on Lemma~\ref{lem:uls-sup}, we obtain the following result on the
performance of the $\operatorname{ULS}$  scan test at supercriticality. As before, we
restrict ourselves to the case where $U_m(t, k_m)$ is of order
$(\log m)^{1/2}$. We also chose to state a simple result instead of a
more precise result with multiple subcases.
This result holds irrespective of the type of boundary condition
assumed on $\Vm$.
\begin{them} \label{thm:uls-sup}
In hypercube detection, the test based on $U_m(t, k_m)$, with $t$
fixed such that $p_c < p_0(t) < 1$ and $\liminf k_m/\log m > 0$ and
$\limsup k_m^{(d-1)/d}/\log m < \alpha d/\delta_{p_0(t)}$, is
asymptotically powerful (resp., powerless) if
\[
\theta_m \bigl[m^{(\alpha-1/2)d} + (\log m)^{d/(2d-2)} \bigr] (\log
m)^{-1/2} \to\infty\qquad\mbox{(resp., $\to0$)}.
\]
\end{them}

We also mention that the equivalent of Theorem~\ref{thm:uls-path}
holds here as well.

The improvement of the supercritical $\operatorname{ULS}$  scan test compared with the
supercritical LOC test is a weaker requirement on $\theta_m$ by a
logarithmic factor. Thus, this test's performance is still much worse
than that of the scan statistic when detecting hypercubes.
%
\subsection{Critical threshold}
If we choose a threshold as described in Section~\ref{sec:cc-cri}, and
if (\ref{S-cri}) is true, then the power of the $\operatorname{ULS}$  scan statistic is greatly
improved, as in the case of the LOC test. In fact, it can be proven
that Theorem~\ref{thm:cc-cube-cri} remains valid with $S(t_m)$
replaced with
$U_m(t_m, k_m)$, as long as $k_m = \mathrm{o}(m)^{\alpha d}$ so that the largest
open cluster under the alternative is scanned. This boils down to
showing that under the null, the $\operatorname{ULS}$  scan statistic is at most a power
of $\log m$,
which we do in Lemma~\ref{lem:uls-orig} below.
However, the $\operatorname{ULS}$  scan test does not seem to offer any substantial gain
in power over the LOC test, given that $\theta_m$ is still required to
be large enough to change the regime of the percolation process within
an alternative $K$ from subcritical to supercritical. That said,
actually proving this would require information on the smaller open
clusters near criticality, which is scarce and very difficult to obtain
(see Borgs \textit{et al.}~\cite{MR1868996} for some partial results
and postulates).
%
\subsection{Across all thresholds}
\label{sec:uls-orig}
Finally, we discuss the (simple) $\operatorname{ULS}$  scan test across all thresholds,
as suggested in Patil and Taillie~\cite{MR2109372}. To take advantage
of a phase
transition near criticality, we assume, as in Section~\ref
{sec:cc-cri}, that
there exists $t_c$ such that $p_0(t_c)=p_c$ and that $p_0(t)$ is a
continuous function of $t$ in a neighborhood of $t_c$. We also assume
that (\ref{S-cri}) holds.
In Proposition~\ref{prp:uls-bad}, we showed that scanning small
clusters may lead to
a decrease in power. For this reason, and also to facilitate the
analysis, we limit ourselves to clusters of size at least $k_m$; that
is, we consider the test based on
%
\begin{equation}\label{uls-orig-mod}
\ULS_m(k_m) = \max_{t \in\bbR} \frac{U_m(t, k_m)}{\sigma_{0|t}},
\end{equation}
where, for definiteness, $U_m(t,k_m)$ is calculated on the torus $\Vm
$ when $t<t_c$.

Let $\Gamma_\theta(\beta) = \inf_{t} \gamma_{\theta|t}(\beta
)/\sigma
_{0|t}^2$, where, in congruence with Sections~\ref{sec:uls-sub} and
\ref{sec:uls-sup},
\[
\gamma_{\theta|t}(\beta) =
\cases{
\gamma\bigl(F_{\theta|t}, \mu_{0|t}, \zeta_{p_\theta(t)}, \beta\bigr), & \quad$t
> t_c$,\cr
\gamma(F_{\theta|t}, \mu_{0|t}, 0, \beta), & \quad$t < t_c$,}
\]
with $\gamma$ being the function defined in Lemma~\ref{lem:gamma}.
%
We first establish the behavior of $\ULS_m(k_m)$ under the null.
%
%
\begin{lem} \label{lem:uls-orig}
Let $k_m = \beta\log m$ where $\beta> 0$, and let $t_\beta$ be such
that $d/\beta\leq\zeta_{p_0(t_\beta)} < \infty$. Define $\eta
(\beta)
:= \sup\{\sigma_{0|t}/\sigma_{0|s} \dvt s \leq t \leq t_\beta\}$. With
probability tending to 1, under $F_0$,
\[
\limsup_{m \to\infty} (\log m)^{-1/2} \ULS_m(k_m) \leq\eta(\beta)
(d \Gamma_{\theta}(\beta))^{1/2}.
\]
If in addition, either $\sigma_{0|t}$ is non-decreasing in $t$ or $F_0$
has no atoms on $(-\infty, t_\beta]$, then, in probability under $F_0$,
\[
(\log m)^{-1/2} \ULS_m(k_m) \to(d \Gamma_{\theta}(\beta))^{1/2}.
\]
\end{lem}

In fact, a result as precise as Lemma~\ref{lem:uls-orig} is
superfluous, given
the behavior of the $\operatorname{ULS}$  scan statistic under the alternative at
supercriticality and near-criticality, which is polynomial in $m$.
The next theorem does not require the use of toroidal boundary conditions.
\begin{them} \label{thm:uls-orig}
In hypercube detection and assuming that (\ref{S-cri}) holds, the test
based on $\ULS_m(k_m)$, with $k_m = [\beta\log m]$ for some $\beta>
0$, is asymptotically powerful if $\theta_m m^\lambda\to\infty$, for
some $0 < \lambda< \alpha/\nu$ satisfying $\lambda<(\alpha-1/2)d$ if
$\alpha> 1/2$.
\end{them}

Thus, scanning all thresholds elicits the best performance of the LOC
tests. Nevertheless, the overall test is still suboptimal when
detecting hypercubes compared with the scan statistic.
We mention in passing that the same result holds for the simpler test
that scans only the largest open cluster at each threshold.
%
\section{Implementation and numerical experiments}
\label{sec:numerics}
The scan test has been shown to be near-optimal in a wide variety of
settings, differing in terms of both network structure and cluster
class (Arias-Castro, Cand\`es and Durand~\cite{cluster}; Arias-Castro,
Donoho and Huo~\cite{MGD}). It is computationally demanding, however.
For the simple situation of detecting a hypercube, the scan statistic
can be computed in $\O(N \log N)$ flops, where $N := m^d$ is the
network size if the size of the hypercube is known. If one scans over
all possible hypercubes, then computing the scan statistic requires $\O
(N^2 \log N)$ flops. For nonparametric shapes, the computational cost
is even higher; in fact, for the problem of detecting a loopless path,
computing the scan statistic corresponds to the reward-budget problem
of DasGupta \textit{et al.} \cite
{DasGuptaHespanhaRiehlSontag06}, shown there to be NP-hard.
Because the scan statistic is so computationally burdensome, the
cluster class is most often taken to be parametric in practice, even
though the underlying clusters may take a much wider range of shapes.
For instance, discs are the prevalent shape used in disease outbreak
detection (Kulldorff and Nagarwalla~\cite{kulldorff1995spatial}), with
variants such as
ellipses (Hobolth, Pedersen and Jensen~\cite{elliptical}; Kulldorff \textit{et al.}~\cite{kulldorff2006ess}). For a wide range of
parametric shapes, Arias-Castro, Donoho and Huo~\cite{MGD} recommended
a multiscale
approximation to the scan statistic.
Efforts to move beyond parametric models include tree-based
approaches (Kulldorff, Fang and Walsh~\cite{tree-based}), simulated
annealing (Duczmal and Assun{\c{c}}{\~a}o~\cite{MR2045632})
and an exhaustive search among arbitrarily shaped clusters of small
size (Tango and Takahashi~\cite{tango2005flexibly}).

The LOC test does not assume any parametric form for the anomalous
cluster, and in that sense is nonparametric. Its computational
complexity at a given threshold is of order the number of nodes plus
the number of edges in the network (Cormen \textit{et al.}~\cite{MR2572804}),
and so of order
$\O(N)$ flops for the square lattice.

The $\operatorname{ULS}$  scan statistic is nonparametric as well. Computing $U
_m(t,k_m)$ requires determining $\Qmt$, which takes $\O(N)$ flops, and
then scanning over $\Qmt$. Because the clusters in $\Qmt$ do not
intersect, scanning over them takes order $\O(N)$ flops. Therefore,
computing $\ULS_m$ can be done in $\O(M \cdot N)$ flops, where $M$ is
the number of distinct values at the nodes. Patil and Taillie
\cite{patil-upper} argued that this can be done faster by using the
tree structure of $\Qm^*$, where the root is the entire network $\Vm$
and a cluster $K \in\Km(t_j)$ is the parent of any cluster $L \in\Km
(t_{j+1})$ such that $L \subset K$, where $t_1 < \cdots< t_M$ denote
the distinct values at the nodes.

We complement our theoretical analysis with some small-scale numerical
experiments. Specifically, we explore the power properties of the LOC
test of Section~\ref{sec:cc} and the $\operatorname{ULS}$  scan test of Section \ref
{sec:uls} in the
context of detecting a hypercube in the two-dimensional square lattice.
Patil, Modarres and Patankar~\cite{uls-soft} are developing sophisticated
software implementing the $\operatorname{ULS}$  scan statistic for use in real-life
situations, with more recent variations Patil, Joshi and Koli \cite
{springerlink:10.1007/s10651-010-0140-1}. However, this software is not
yet available, so we implemented our own (basic) routines.

We used the statistical software \textit{R} (R Core Team~\cite{R}) with
the package \texttt{igraph} (Csardi~\cite{igraph}). Our (basic)
implementation of the $\operatorname{ULS}$  scan statistic for a given threshold is much
slower than both the scan statistic with a given mask and the LOC
statistic, especially when there is no constraint on the size of the
open clusters to be scanned, that is, when $k_m = 1$. In all of our
experiments, we chose the square lattice in dimension $d=2$ with side
length $m = 500$ for a total of 250,000 nodes, and we considered three
alternatives: squares of side length $\ell\in\{10, 50, 100\}$,
corresponding roughly to $\alpha\in\{0.4, 0.7, 0.8\}$. The squares
were fixed away from the boundary of the lattice, given that the
methods are essentially location-independent. (This is rigorously true
of the scan statistic.) We assessed the performance of a method in a
given situation by estimating its risk, which we define as the sum of
the probabilities of type I and type II errors optimized over all
rejection regions.
%

\begin{figure}

\includegraphics{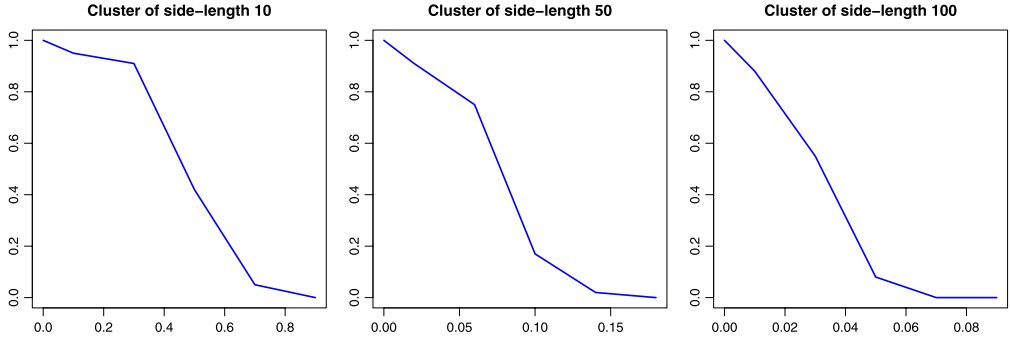}

\caption{The risk of the scan test against each of the three
alternatives. The $x$-axis is $\theta$, and the $y$-axis is the
estimated risk based on 100 replicates.}\label{fig:scan}
\end{figure}

We first ran some experiments to quickly assess the power of the scan
test. We found that the test agrees very well with the theory (i.e.,
Lemma~\ref{lem:scan-cube}), which we already knew from previous experience.
Specifically, we assumed a normal location model and simulated 100
realizations of the null and each of the three alternatives with
$\theta\in\{j/\ell\dvt j=1,3,5,7,9\}$ (see Figure~\ref{fig:scan}).
%
%
\begin{figure}[b]

\includegraphics{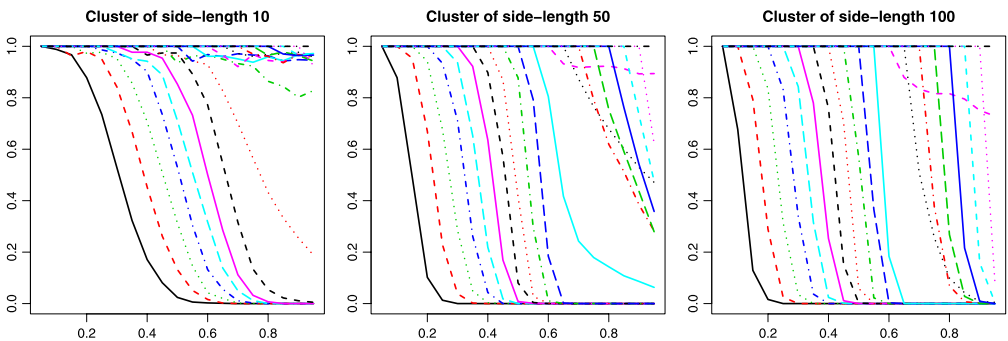}

\caption{The risk of the LOC test against each of the three
alternatives. The $x$-axis is the percolation probability $q$ on the
anomalous cluster, and the $y$-axis is the estimated risk based on
1000 replicates. Each curve corresponds to a different percolation
probability $p$.}\label{fig:cc}
\end{figure}

\begin{figure}

\includegraphics{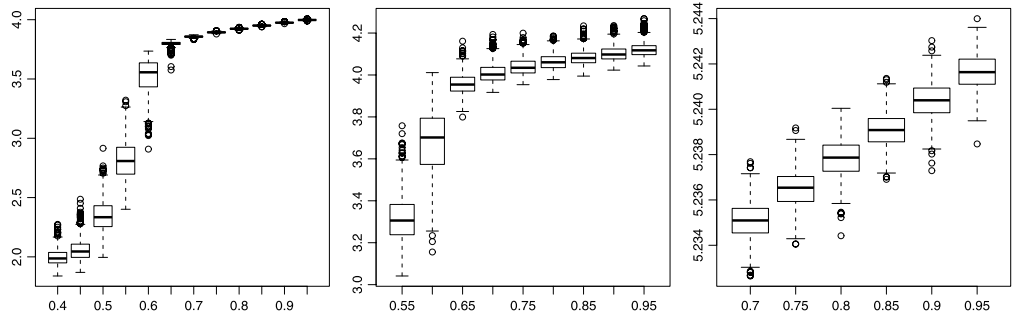}

\caption{The size of the largest open cluster in $\log_{10}$ scale
($y$-axis) versus the percolation probability $q$, for the alternative
$\ell= 100$ and $p \in\{0.40, 0.55, 0.70\}$ (from left to right).
Each boxplot represent 1000 replicates.}\label{fig:cc-box}
\end{figure}

Next, we performed some larger experiments to assess the power of the
LOC test. We simply assumed a site percolation model with probability
$p \in\{0.05, 0.10, \ldots, 0.90, 0.95\}$. Note that $p_c$ is not known
for site percolation in the square lattice, although $p_c \approx
0.593$ from extensive numerical experiments (Feng, Deng and Bl\"ote
\cite{PhysRevE.78.031136}).
We simulated the null and each of the three alternatives with $q \in\{
0.05, 0.10, \ldots, 0.90, 0.95\}, q > p,$ within the anomalous
cluster. We replicated each situation 1000 times. The risk curves are
shown in Figure~\ref{fig:cc}. The test seems to behave similarly above and
below criticality. At near-criticality, the test is rather erratic.
However, when the size of the anomalous cluster is large enough, $\ell
= 100$, the risk curve is steepest just under $p_c$, at $p =0.55$ in
our experiments, with full power against $q \geq0.65$. Figure~\ref{fig:cc-box}
shows boxplots of the test statistic for the case where $\ell= 100$
and $p = 0.40$ (subcritical), $p = 0.55$ (near-critical), and $p =
0.70$ (supercritical).
%

\begin{figure}[b]

\includegraphics{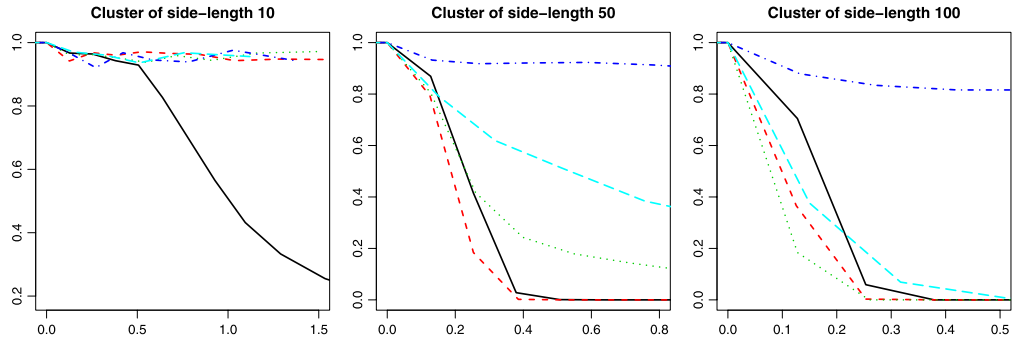}

\caption{The risk of the LOC test in the context of a normal location
model. The $x$-axis is $\theta$, and the $y$-axis is the estimated risk
based on 1000 replicates. Each curve corresponds to a different
threshold~$t$. The solid (---), dashed (-\,-), dotted ($\cdots$), dot-dashed (-$\cdot$-) and long-dashed (--\,--) curves correspond
to $p = 0.40, 0.50, 0.55, 0.60 \mbox{ and } 0.70$, respectively.}\label{fig:cc-norm}
\end{figure}

If we were to use this test in the context of a normal location model,
then the correspondence would be $t = \bar{\Phi}^{-1}(p)$ (the
threshold) and $\theta= t -\bar{\Phi}^{-1}(q)$, where $\bar{\Phi}$
denotes the normal survival distribution function. Figure~\ref{fig:cc-norm}
plots the risk curves in this context for $p \in\{0.40, 0.50, 0.55,
0.60, 0.70\}$. In particular, the test at near-criticality with $t =
\bar{\Phi}^{-1}(0.55) = -0.126$ has full power against the alternative
with $\ell= 100$ and $\theta= 0.26$.
%

%
%
%
%
%

Finally, we experimented with the $\operatorname{ULS}$  scan test. To limit the size of
our simulations, we considered alternatives with $\theta= \Phi
^{-1}(q)$ with $q \in\{0.55, 0.6, 0.65, 0.70, 0.80, 0.90\}$ and chose
$t = \Phi^{-1}(p)$ with $p \in\{0.40, 0.50, 0.55, 0.60, 0.70\}$ as
thresholds. We restricted scanning to open clusters of size not smaller
than $1/10$ of the size of largest open cluster, essentially falling in
the regime of Part 2 of Lemma~\ref{lem:uls-sub}, and also making the
computation much faster. We used 200 replicates. We again see that the
risk curve is sharpest near criticality when the size of the anomalous
cluster is sufficiently large, here for $\ell\geq50$. Compared with
the LOC test, the $\operatorname{ULS}$  scan test has more power at large $\theta$ when
the cluster is small $\ell= 10$ (as predicted) and, more
interestingly, slightly more power when the cluster is larger. Compared
with the scan statistic, which knows the size and shape of the anomalous
cluster, the $\operatorname{ULS}$  scan test with the best choice of threshold
(corresponding to $p = 0.55$) requires approximately threefold greater
signal amplitude.
%
%
%
\begin{figure}

\includegraphics{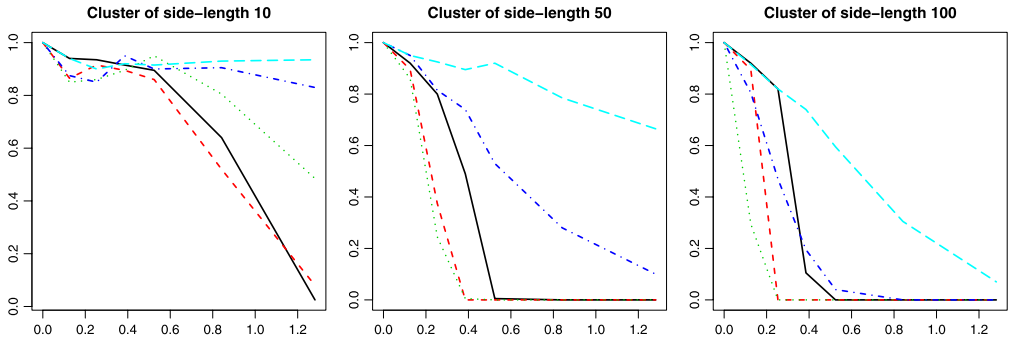}

\caption{The risk of the $\operatorname{ULS}$  scan test against each of the three
alternatives. On the $x$-axis is $\theta$, and on the $y$-axis is the
estimated risk based on 200 replicates. Each curve corresponds to a
different threshold~$t$. The solid (---), dashed (-\,-), dotted ($\cdots$), dot-dashed (-$\cdot$-) and long-dashed (--\,--) curves
correspond to $p = 0.40, 0.50, 0.55, 0.60 \mbox{ and } 0.70$,
respectively.}\label{fig:uls}
\end{figure}
%
\section{Discussion}
\label{sec:discussion}
The contribution of this paper is a rigorous mathematical analysis of
the performance of the LOC test independent of, and more extensively
than Davies, Langovoy and
Wittich~\cite{langovoy-10} and Langovoy and Wittich~\cite{langovoy-11}, and of the $\operatorname{ULS}$  scan test, both
nonparametric and computationally tractable methods. We made abundant
use of percolation theory to establish these results. We compared the
power of these tests with that of the scan statistic, which is known to
be near-optimal in a wide array of settings. Although these tests are
comparable in power with the scan statistic for the detection of a
path, they may be substantially less powerful for the detection of a
hypercube. Note, however, that the scan statistic is provided with
knowledge about the shape and size of the anomalous cluster. In theory,
we argued that this was the case based on some heuristics and
conjectures from percolation theory. Numerically, this appears to be
the case when the anomalous cluster is large enough. In our
experiments, the $\operatorname{ULS}$  scan test was slightly more powerful than the LOC
test, and required a $\theta$ three to four times larger than the scan
statistic, which has the advantage of knowing the shape and size of the
cluster. This result is promising, and further numerical experiments
are needed to evaluate the power of these tests in truly nonparametric
settings, because they do not require previous information about
cluster shape, and are computationally more feasible in general.

Our theoretical results generalize to other networks that resemble the
lattice, with a different critical percolation probability $p_c$ and
different functions $\zeta_p$ and $\delta_p$. In particular, we used
the self-similarity property of the square lattice and the fact that it
has polynomial growth. Our results also generalize to other cluster
classes; in the setting of the square lattice, they extend immediately
to any class of clusters that includes a hypercube of comparable size
(e.g., the class\vadjust{\goodbreak} $\Km$ of clusters $K$ of size $|K| = [m^\alpha]^d$),
such that there is a hypercube $K_0 \subset K$ with $|K_0|/|K| \geq
\omega_m$, where $\omega_m \to0$ more slowly than any negative power
of $m$. In addition, the class might contain clusters of different
sizes, although in that case the worst-case risk would be driven by the
smallest clusters. Implementation of the scan statistic may be much
more demanding in this case.
The main results of Section~\ref{sec:cc} require only that $F_\theta
(t)$ be
twice differentiable in $(t, \theta)$, with $\partial_\theta F_\theta
(t) < 0$ for all $(t, \theta)$, which is the case, for example, for
location models and scale models if $F_0$ is twice differentiable with
a strictly positive first derivate.
With some additional work, we also can obtain results for classes of
``thin'' clusters as defined in Arias-Castro, Cand\`es and Durand \cite
{cluster}. The key is to understand
the percolation behavior within and near such clusters. Some results
are available for slabs (Grimmett~\cite{MR1707339}, Theorem 7.2)
and more general subgraphs of lattices including ``wedges,'' and these
appear to be transferable to other ``curved'' slabs.
%
%
%
%
\begin{appendix}\label{appa}
\section*{Appendix A: Proofs}
\setcounter{equation}{0}
\renewcommand{\thelem}{A.\arabic{lem}}
\begin{nota*}
We write $f_m \sim g_m$ as $n\to\infty$ if $f_m/g_m\to1$. Similarly,
we use $\O(\cdot)$ and $\o(\cdot)$ and write $f_m \asymp g_m$ as $n\to\infty$ if
$f_m = \O(g_m)$ and vice versa. We also use their random counterparts,
$\simP$, $\asymP$, $\OP(\cdot)$, and $\o_\mathrm{P}(\cdot)$.
For example, $Z_m=\o_\mathrm{P}(k_m)$ means that $Z_m/k_m\to0$ in
probability,
and $Z_m=\OP(k_m)$ means that $Z_m/k_m$ is bounded in probability,
which is to say that $\P(|Z_m| \ge k_m l_m) \to1$ as $m \to\infty$
for any $l_m$ satisfying $l_m \to\infty$.
We use $1\{A\}$ to denote the indicator function of the set $A$.
The maximum of $k$ and $\ell$ is denoted by $k \vee\ell$.
\end{nota*}
\subsection{On the size of percolation clusters}
Here we state and prove some results on the sizes of percolation
clusters in $\bbZ^d$.
We start by proving some properties of $\zeta_p$.
Recall that $S$ denotes the size of the open cluster at the origin.
Besides the limit in (\ref{zeta}), the following bound holds for
$p<p_c$ and all $k \geq1$:
%
\begin{equation}\label{S-ineq}
\P_p(S \geq k) \leq(1-p)^2 \frac{k \mathrm{e}^{-k \zeta_p}}{(1 - \mathrm{e}^{-\zeta_p})^2},
\end{equation}
by Grimmett~\cite{MR1707339}, Equation (6.80), adapted to site percolation.
\setcounter{lem}{0}
\begin{lem} \label{lem:zeta}
The function $\zeta_p$ defined in (\ref{zeta}) is continuous and
strictly decreasing over $(0, p_c]$, and satisfies $\lim_{p \to0}
\zeta
_p = \infty$ and $\lim_{p \to p_c} \zeta_p = 0$.
\end{lem}
\begin{pf}
Let $0 \le p < p' \le1$. By coupling $\bbP_p$ and $\bbP_{p'}$ in the
usual way,
\[
\bbP_p(S=k) \ge(p/p')^k \bbP_{p'}(S=k),
\]
so that $\zeta_p \le\zeta_{p'} + \log(p'/p)$.
Applying Grimmett~\cite{MR1707339}, Theorem 2.38, to the event $\{S\ge k\}
$, we find
that, as in the proof of Grimmett~\cite{MR1707339},\vadjust{\goodbreak} Equation (6.16), $\zeta
_p/\log p
\le\zeta_{p'}/\log p'$. In summary,
%
\begin{equation}\label{G2}
\zeta_p \biggl(1-\frac{\log(1/p')}{\log(1/p)} \biggr) \le\zeta_p - \zeta
_{p'}\le\log(p'/p).
\end{equation}
Therefore, $\zeta_p$ is continuous and strictly decreasing on $(0,p_c)$.
Moreover, by fixing $p' \in(0, p_c)$ and letting $p \to0$, we have
\[
\zeta_p \geq\zeta_{p'} \frac{\log(1/p)}{\log(1/p')} \to\infty.
\]

Finally, by Grimmett~\cite{MR1707339}, Equations (6.83), (6.56), $\zeta_p
\to
0 =
\zeta_{p_c}$ as $p \uparrow p_c$.
\end{pf}

Next, we prove (\ref{S-conv}). We do this by standard means, and
the claim may be strengthened (see also Grimmett~\cite{MR886037};
Hofstad and Redig~\cite{vdHR}).
\begin{lem} \label{lem:S-conv}
Consider site percolation on $\bbZ^d$ with parameter $p < p_c$, and let
$S_m$ denote the size of the largest open cluster within $\Vm$. Then
(\ref{S-conv}) holds, namely
\[
\frac{S_m}{\log m} \to\frac{d}{\zeta_p}, \qquad\mbox{in probability}.
\]
\end{lem}
\begin{pf}
Fix $0 < \eps< 1/2$. Let $S^{v}$ be the size of the open cluster at a
node $v \in\bbZ^d$, which has the same distribution as $S$.
We start with the upper bound. By the union bound,
%
\begin{equation}\label{Sm-ub}
\mathbb{P}(S_m \geq k) \leq\sum_{v \in\Vm} \mathbb{P}(S^v \geq k)
= |\Vm| \cdot\mathbb{P}(S \geq k).
\end{equation}
Thus, using (\ref{zeta}), for $k_m(\eps) := (1 +\eps) (d/\zeta_p)
\log
m$ and $m$ large enough,
\[
\mathbb{P}\bigl(S_m \geq k_m(\eps)\bigr) \leq m^d \exp\bigl(- (1-\eps/2)\zeta_p
k_m(\eps)\bigr)
\leq m^{-\eps d/4},
\]
and the term on the right-hand side converges to $0$.

For the lower bound, consider $N = \lceil m^d/(\log m)^{2d}\rceil$
nodes $v_1, \ldots, v_N \in\Vm$ separated from each other and the
boundary of $\Vm$ by at least $\frac12 (\log m)^2$.
Let $k_m(\eps) := (1 -\eps) (d/\zeta_p) \log m$. For sufficiently large
$m$, the events $E_i : = \{|S^{v_i}|\le k_m(\eps)\}$
are independent.
Therefore, using (\ref{zeta}), for large $m$,
%
\begin{eqnarray}\label{Sm-lb}
\mathbb{P}\bigl(S_m \leq k_m(\eps)\bigr)
&\leq&\bigl(1 -\mathbb{P}\bigl(S \geq k_m(\eps)\bigr)\bigr)^N  \nonumber\\
&\leq&\bigl(1 -\exp\bigl(- (1+\eps/2)\zeta_p k_m(\eps)\bigr)\bigr)^N \\
&\leq&\exp\bigl(-m^{\eps d/2}/(\log m)^{2d}\bigr), \nonumber
\end{eqnarray}
and the last term on the right-hand side tends to 0 as $m \to\infty$.
\end{pf}

The following result describes the behavior of size of the open cluster
at the origin when $p$ is small.
It may be made more precise, but we do not pursue this here.
\begin{lem} \label{lem:p-small}
There exists $c > 0$ depending only on $d$ such that, for $p \in(0,
(2c)^{-1})$,
\[
p^k \le\bbP_p(S \geq k) \le\tfrac12 (c p)^k\qquad \forall k \geq1.
\]
%
\end{lem}
\begin{pf}
An \emph{animal} is a connected subgraph of $\bbZ^d$ containing the origin.
The lower bound comes from considering the probability that any given
animal of size $k$ is open.
For the upper bound, by the union bound, we have $\bbP_p(S= k) \le
|\cA
_k| p^k$, where $\cA_k$ is the set of animals with $k$ vertices. There
is a constant $c > 0$ such that $|\cA_k| \le c^k$, so that
\[
\bbP_p(S \geq k) \leq\sum_{\ell\geq k} c^\ell p^{\ell} = \dfrac{(c
p)^k}{1 - cp}
\le\frac12 (cp)^k,
\]
when $cp<\frac12$.
\end{pf}

We next present a result on the number of open clusters of a given size
that is valid for all $p \in(0,1)$.
%
\begin{lem} \label{lem:Nk}
Consider site percolation on $\bbZ^d$ with parameter $p$, and let
$N_m(k)$ denote the number of open clusters of size $k$ within $\Vm$.
Then, for $k \ge1$,
\[
\frac{(m -2k)^d}{k} \mathbb{P}(S=k) \leq\mathbb{E}(N_m(k)) \leq
\frac{m^d}k \mathbb{P}(\infty> S\geq k),
\]
In addition, for $k,\ell\ge1$,
\[
\bigl|\operatorname{Cov}(N_m(k), N_m(\ell))\bigr| \leq
3^{d+1}(k+\ell)^d \mathbb{E}(N_m(k \vee\ell)).
\]
Thus, for $k \ge1$,
\[
\operatorname{Var}(N_m(k)) \leq6^{d+1}k^d \mathbb{E}(N_m(k)).
\]
\end{lem}
\begin{pf}
%
Let $S^v_m$ be the size of the open cluster at $v$ within the box $\Vm
$. Then
%
\begin{equation}\label{G3}
N_m(k) = \sum_{v \in\Vm} X^v(k),
\end{equation}
where $X^v(k) = k^{-1}1\{S^v_m=k\}$.
We immediately have
\[
\mathbb{E}(N_m(k)) \leq\sum_{v \in\Vm} \frac1k \mathbb{P}(\infty
>S^v \geq k) =
\frac{|\Vm|}k \mathbb{P}(\infty>S \geq k).
\]
For the lower bound, we count only nodes away from the boundary, obtaining
\[
\mathbb{E}(N_m(k)) \geq|\Vm(k)| \frac1k \mathbb{P}(S=k),
\]
where $\Vm(k):= \{k,\ldots,m-k\}^d$.\vadjust{\goodbreak}

We turn now to the covariances. By (\ref{G3}),
\begin{eqnarray*}
\operatorname{Cov}(N_m(k),N_m(\ell)) &=& \sum_{v,w\in\Vm}
\operatorname{Cov}(X^v(k),X^w(l))\\
&=& \mathop{\mathop{\sum}_{v,w\in\Vm}}_{\|v-w\|\le
k+\ell} \operatorname{Cov}(X^v(k),X^w(l)),
\end{eqnarray*}
because $X^v(k)$ and $X^w(\ell)$ are independent if $\|v-w\|> k+\ell$,
where $\|\cdot\|$ denotes $\ell^\infty$-norm.
Now,
\begin{eqnarray*}
\bigl|\operatorname{Cov}(X^v(k),X^w(\ell))\bigr| &=&\bigl|\mathbb
{E}\bigl(X^w(\ell)\vert X^v(k) = k^{-1}\bigr) - \mathbb{E}(X^w(\ell))\bigr|
\mathbb{E}(X^v(k))\\
&\le&\frac1\ell\mathbb{E}(X^v(k)),
\end{eqnarray*}
so that
\[
\bigl|\operatorname{Cov}(N_m(k),N_m(\ell))\bigr| \le\frac1\ell
(2k+2\ell+1)^d
\mathbb{E}(N_m(k)),
\]
and the second claim of the lemma follows.
\end{pf}

We now describe some properties of the open clusters within $\Vm$ in
the supercritical regime.
In this regime, it is known that, with probability 1, there is a unique
infinite open cluster in $\bbZ^d$, denoted by $Q_\infty$ (see, e.g.,
Grimmett~\cite{MR1707339}, Section 8.2). With high probability, the
largest open
cluster within $\Vm$ is a subgraph of this infinite open cluster. Next,
we present some additional information on its size, $S_m$.
\begin{lem} \label{lem:large-sup}
Suppose that $p > p_c$. There is a constant $C>0$ such that, with
probability at least $1 -\exp(-C m^{d-1})$, there is a unique largest
open cluster within $\Vm$, and it is a subgraph of $Q_\infty$.
Moreover, as $m \to\infty$, its size $S_m$ satisfies
\[
\frac{S_m -\mathbb{E}(S_m)}{\sqrt{\operatorname{Var}(S_m)}} \to\cN
(0,1), \qquad\mbox{in
distribution},
\]
with $\mathbb{E}(S_m) \sim\Theta_p |\Vm|$ and $\operatorname
{Var}(S_m) \sim\sigma^2
|\Vm|$ for some $\sigma^2 > 0$ depending on $(d,p)$.
\end{lem}
\begin{pf}
For the first part and the limiting behavior of $\mathbb{E}(S_m)$ as $m
\to\infty$, see the discussion
of Penrose and Pisztora~\cite{MR1372330}, Theorems 4 and 6, and the
beginning of this Appendix.
For the weak limit and the limit size of the variance of $S_m$, see,
for example, Penrose~\cite{MR1880230}, Theorem 3.2.
\end{pf}

We next describe some properties of the smaller open clusters. Let
$\Stm
$ be the size of the largest open cluster
of $\bbZ^d$ that is contained entirely within $\Vm$.
\begin{lem} \label{lem:small-sup}
Suppose that $p > p_c$. There exists a positive constant $\delta_p$
such that
\[
\frac{\Stm}{(\log m)^{d/(d-1)}} \to\biggl(\frac d{\delta_p}
\biggr)^{d/(d-1)}, \qquad\mbox{in probability}.
\]
For any $c >0$, there exists $\sigma_i=\sigma_i(p,c)>0$ such that the
following holds:
With probability tending to 1, there exist at least
$\sigma_1 m^d \exp[-\sigma_2(\log m)^{(d-1)/d}]$ open clusters
of size $[c \log m]$ of $\bbZ^d$ lying within $\Vm$.
\end{lem}

Our results on exact asymptotics in the supercritical phase concern
$\Vm
$ with toroidal boundary conditions. One effect of removing the
boundary from $\Vm$ is that the asymptotics of the largest cluster
coincide with those of $S_m$, as well as for the second-largest cluster
$S_m^{(2)}$. In the proof of Theorem~\ref{thm:uls-sup}, we need an
upper bound on the size of the second-largest cluster inside a box with
``free'' boundary conditions. We do not explore this in detail here,
because it relies on extensions of arguments of Kesten and Zhang \cite
{kestenz} (see
also Grimmett~\cite{MR1707339}, Proof of Theorem 8.65), which have not yet
been not
fully explored in the literature. Instead, we note that the the
second-largest open cluster in a supercritical percolation model on
$\Vm
$ with free boundary conditions has size of order $\OP((\log m)^{d/(d-1)})$.
\begin{pf*}{Proof of Lemma~\ref{lem:small-sup}}
It was proven by Cerf~\cite{Cerf06} that the limit
%
\begin{equation}\label{delta2}
\delta_p : = -\lim_{k\to\infty} k^{-(d-1)/d} \log\bbP(S=k)
\end{equation}
exists and is strictly positive and finite when $p_c < p < 1$. It is elementary
that $\delta_p$ thus defined is equal to that of
(\ref{delta}) (see also Grimmett~\cite{MR1707339}, Section 8.6).
The first part of the lemma follows by the same proof as used in Lemma
\ref{lem:S-conv}.

As in the proof of Lemma~\ref{lem:Nk}, the mean number $\mu_m$ of
clusters of
size $k:=[c\log m]$ satisfies
\[
\frac{m^d}{c\log m} \exp\bigl(-\delta^1(c\log m)^{(d-1)/d} \bigr) \le
\mu_m
\le\frac{m^d}{[c\log m]} \exp\bigl(-\delta^2(c\log m)^{(d-1)/d} \bigr)
\]
for positive constants $\delta^i$. The number of such clusters has
variance no larger than
$Ck^d\mu_m$ for some $C<\infty$. The claim follows by Chebyshev's inequality.
\end{pf*}
%
%
%
%
%
\subsection{Some distributional properties}
Here we present some results for $\operatorname{AEP}$ and exponential
families of
distributions. Our first result is on the size of the maximum of an
i.i.d. sample from an $\operatorname{AEP}$ distribution.
\begin{lem} \label{lem:aep-max}
Let $F \in\operatorname{AEP}(b, C)$ for some $b > 0$ and $C > 0$.
Then, for
$X_1, \ldots, X_n \stackrel{\mathrm{i.i.d.}}{\sim} F$,
\[
\frac{\max(X_1, \ldots, X_n)}{(\log n)^{1/b}} \to C^{-1/b}, \qquad\mbox
{in probability}.
\]
\end{lem}
\begin{pf}
Fix $\eps\in(0,1)$ and define $x_n(\eps) = ((1-\eps) (\log
n)/C)^{1/b}$. For $n$ large enough, we have, by independence,
\begin{eqnarray*}
\mathbb{P}\bigl(\max(X_1, \ldots, X_n) \leq x_n(\eps)\bigr)
&\leq&\bigl(1 -\bar{F}(x_n(\eps))\bigr)^{n} \\
&\leq&\bigl(1 -\exp\bigl(-(1+\eps) C x_n(\eps)^b\bigr)\bigr)^{n} \\
&\leq&\exp(-n^{\eps^2}) \to0.
\end{eqnarray*}
Now redefine $x_n(\eps) = ((1+\eps) (\log n)/C)^{1/b}$.
For $n$ large enough, we have, by the union bound,
\begin{eqnarray*}
\mathbb{P}\bigl(\max(X_1, \ldots, X_n) \geq x_n(\eps)\bigr)
&\leq& n \bar{F}(x_n(\eps)) \\
&\leq& n \exp\bigl(-(1-\eps/3) C x_n(\eps)^b\bigr) \\
&\leq& n^{-\eps/3} \to0.
\end{eqnarray*}
\upqed\end{pf}

We next describe the behavior at infinity of the logarithmic
moment-generating function and rate function of an $\operatorname{AEP}$
distribution.
\begin{lem} \label{lem:aep-rate}
Let $F \in\operatorname{AEP}(b, C)$ for some $b \geq1$ and $C > 0$, with
logarithmic moment-generating function $\Lambda$ and rate function
$\Lambda^*$. Then, as $\theta\to\infty$,
%
\begin{eqnarray}
\theta^{-b/(b-1)} \Lambda(\theta) &\to&C (b-1) (C b)^{-b/(b-1)}, \qquad b >1; \label{Lda} \\
\bigl(\log\bigl(1/(C-\theta)\bigr)\bigr)^{-1} \Lambda(\theta) &\to&1,\qquad b = 1;
\label{Lda1}
\end{eqnarray}
and,  as $x \to\infty$,
%
\begin{equation}\label{Lda-star}
x^{-b} \Lambda^*(x) \to C.
\end{equation}
\end{lem}
\begin{pf}
Let $\varphi$ be the moment-generating function of $F$. We focus on the
upper bound in (\ref{Lda}) -- obtaining the bound in (\ref{Lda1}) is
analogous -- and deduce the lower bound in (\ref{Lda-star}). Let $b>1$,
$C/2 < A < C$, and let $x_1 > 0$ be such that $\bar{F}(x) \leq\exp(-A
x^b)$ for all $x > x_1$. We start from the following bound:
\[
\varphi(\theta) = \int_{-\infty}^\infty\theta\exp(\theta x) \bar
{F}(x) \,\mathrm{d}x \leq\exp(\theta x_1) + \int_{x_1}^\infty\theta\exp
(\theta
x -A x^b)\,\mathrm{d}x.
\]
We again divide the integral into $x \leq x_2$ and $x > x_2$, where
$x_2 :=(2 \theta/A)^{1/(b-1)}$. For $x \leq x_2$, we bound
$\exp
(\theta x -A x^b)$ by its maximum over $(0, \infty)$. For $x > x_2$,
$\exp(\theta x -A x^b) \leq\exp(-(C/4) x^b)$. Letting $B = A (b-1) (A
b)^{-b/(b-1)}$, and assuming that $\theta$ is large enough such
that $x_2 > x_1$, we get
\[
\int_{x_1}^\infty\theta\exp(\theta x -A x^b)\,\mathrm{d}x
\leq(x_2 -x_1) \theta\exp\bigl( B \theta^{b/(b-1)} \bigr) +
\theta\int_{x_2}^\infty\exp\bigl(-(C/4) x^b\bigr)\,\mathrm{d}x.
\]
Thus, when $\theta\to\infty$,
%
\begin{equation}\label{mgf-ub}
\varphi(\theta) = \O\bigl(\theta^{b/(b-1)}\bigr) \exp\bigl( B \theta
^{b/(b-1)} \bigr).
\end{equation}
Taking logs and letting $\theta\to\infty$, we get
\[
\limsup_{\theta\to\infty} \theta^{-b/(b-1)} \Lambda(\theta)
\leq A (b-1) (A b)^{-b/(b-1)}.
\]
Then letting $A$ tend to $C$, we obtain the upper bound in (\ref{Lda}).

Now, for $x$ exceeding the mean of $F$, $\Lambda^*(x) = \sup_{\theta
\geq0} (\theta x - \Lambda(\theta))$, and starting from (\ref
{mgf-ub}), we obtain
\[
\Lambda^*(x) \geq\sup_{\theta\geq0} \bigl(\theta x -B \theta^{b/(b-1)}\bigr)
-\log2 = A x^b -\log2.
\]
Therefore,
\[
\mathop{\underline{\lim}}_{x\to\infty} x^{-b} \Lambda^*(x) \geq A.
\]
Then, letting $A$ tend to $C$, we obtain the lower bound in (\ref{Lda-star}).
\end{pf}
%
%

We now define $\gamma$, first appearing in Section~\ref{sec:uls-sub}.
Our function $\gamma$ depends on certain quantities listed in the
following lemma.
It also depends on the quantity $\zeta$, which we take as that defined
in (\ref{zeta}).
It is only through its dependence on $\zeta$ that $\gamma$ is affected
by the geometry of $\Vm$.
%
%
\begin{lem} \label{lem:gamma}
Consider a distribution $F$ on the real line, possibly discrete but not
a point mass, with finite mean $\mu$ and finite moment-generating
function at some positive $\theta> 0$, and let $\Lambda^*$ denote its
rate function. Let $\nu\leq\mu$, and fix $\beta, \zeta\in
[0,\infty)$.
\begin{enumerate}
\item Assume that $\zeta\ne0$. If $0 < \beta< 1/\zeta$, or $\beta=
0$ and $F \in\operatorname{AEP}(b, C)$ for some $b \geq2$ and $C >
0$, then
there is a unique solution $\gamma=\gamma(F, \nu, \zeta, \beta)$
to the
following equation
\[
\inf_{\beta< s < 1/\zeta}\bigl [s \Lambda^* \bigl(\nu+ \sqrt{\gamma
/s} \bigr) + s \zeta\bigr] = 1.
\]
\item Assume that $\zeta=0$. The foregoing holds as long as $\nu=\mu$
(and with $1/\zeta$ interpreted as $\infty$).
\end{enumerate}
\end{lem}
\begin{pf}
Let $M=\sup\{x\dvt \Lambda^*(x)<\infty\}$. Because $F$ is not a point
mass, $\mu< M \le\infty$.
Define
\[
G(s, \gamma) = s \Lambda^* \bigl(\nu+ \sqrt{\gamma/s} \bigr) + s \zeta.
\]
Note that $G(s,\gamma)$ is finite (resp., infinite) if
$\gamma/s < (M-\nu)^2$ (resp., $\gamma/s > (M-\nu)^2$).
In addition, $G(s,\gamma)$, and its derivatives are continuous wherever
$G$ is finite, and thus are uniformly continuous
on any compact subset of $[0,\infty)^2$ on which $G$ is finite. Furthermore,
$G(s,\gamma)$ is strictly increasing in $\gamma$ on the interval
$(0,s(M-\nu)^2)$.
Let
%
\begin{equation}\label{v4}
L_\beta(\gamma) = \inf_{\beta< s < 1/\zeta} G(s, \gamma).
\end{equation}

Thus $L_\beta(\gamma)$ is finite if $\gamma\zeta< (M-\nu)^2$, and
infinite when
$<$ is replaced by $>$. Furthermore, for $\gamma<(M-\nu)^2/\zeta$, the
infimum is achieved at some value $s_\gamma$ of $s$ in a neighborhood
where $G(s,\gamma)<\infty$.

Assume first that $\beta> 0$. It may be seen that $L_\beta(\gamma)$ is
continuous and strictly increasing in $\gamma$ on the interval
$[0,(M-\nu)^2/\zeta)$. Let $0\le\gamma<\gamma'<(M-\nu)^2/\zeta$. Then
%
\begin{equation}\label{L-ub}
0 \le L_\beta(\gamma') - L_\beta(\gamma) \le G(s_\gamma,\gamma
')-G(s_\gamma,\gamma),
\end{equation}
and continuity follows from the properties of $G$ noted earlier. Similarly,
%
\begin{equation}\label{L-lb}
L_\beta(\gamma') -L_\beta(\gamma) \geq G(s_{\gamma'}, \gamma') -
G(s_{\gamma'},\gamma)
\end{equation}
and strict monotonicity follows similarly.

It suffices to prove that $L_\beta(\gamma)$ takes values $<$1 and \emph
{finite} values $>$1. The first claim follows from the fact that,
with $\gamma= \beta(\mu-\nu)^2$,
\[
L_\beta(\gamma) \leq G(\beta, \gamma) = \beta\zeta< 1.
\]

We now turn to the second claim, and make use of two general properties
of rate functions that follow from Dembo and Zeitouni~\cite{MR2571413},
Equation (2.2.10), Lemma 2.2.20.
It is standard that $\Lambda^*(\mu+x) \sim\frac12 (x/\sigma)^2 $
as $x
\downarrow0$, where
$\sigma^2>0$ is the variance of $F$. Therefore,
%
\begin{equation}\label{v3}
\exists T\in(0,M) \mbox{ such that } \Lambda^*(\mu+x) \ge\tfrac14
(x/\sigma)^2
\mbox{ when } 0 \le x \le T.
\end{equation}
With $T$ thus chosen, by convexity,
%
\begin{equation}\label{v2}
\exists A>0 \mbox{ such that } \Lambda^*(\mu+x) \ge A x
\mbox{ when } x \ge T.
\end{equation}

Assume first that $\zeta>0$ and $M=\infty$.
By (\ref{v2}), for sufficiently large $\gamma$,
\[
\infty> L_\beta(\gamma) \geq\inf_{\beta< s < 1/\zeta} \bigl[sA  \bigl(\nu
-\mu+ \sqrt{\gamma/s} \bigr) + s \zeta \bigr] \ge A\bigl(\beta(\nu-\mu
)+\sqrt
{\gamma\beta}\bigr) > 1.
\]

Suppose next that $\zeta>0$ and $M<\infty$.
Let $0 < \gamma< (M-\nu)^2/\zeta$. Because $\Lambda^*(\nu+\sqrt
{\gamma
/s})=\infty$ if $s<\gamma/(M-\nu)^2=:\beta_0(\gamma)$,
%
\begin{eqnarray}\label{v1}
\infty&>& L_\beta(\gamma) \ge\beta_0\inf_{\beta_0<s<1/\zeta}
\Lambda
^* \bigl(\nu+ \sqrt{\gamma/s} \bigr) +\beta_0\zeta\nonumber
\\[-8pt]
\\[-8pt]
&=& \beta_0 \Lambda^*\bigl(\nu+ \sqrt{\gamma\zeta}\bigr) +\beta_0\zeta.\nonumber
\end{eqnarray}
The limit of this, as $\gamma\uparrow(M-\nu)^2/\zeta$, is strictly
greater than $1$.

Now let $\zeta=0$ and $\nu=\mu$, and note that $L_\beta(\gamma
)<\infty$
for all $\gamma\ge0$. Suppose that $M\le\infty$ and
$\gamma>0$.
By dividing the infimum in (\ref{v4}) according to whether or not
$\sqrt
{\gamma/s}<T$, we find that
\begin{eqnarray*}
\infty&>& L_\beta(\gamma) \ge
\min\Bigl\{\inf_{\beta<s<\gamma/T^2} s\Lambda^*\bigl(\mu+\sqrt{\gamma/s}\bigr),
\inf_{s>\gamma/T^2} s\Lambda^*\bigl(\mu+\sqrt{\gamma/s}\bigr) \Bigr\}\\
&\ge&\min\biggl\{ A\sqrt{\gamma\beta}, \frac14 \gamma/\sigma^2 \biggr\},
\end{eqnarray*}
by (\ref{v3})--(\ref{v2}).
This diverges as $\gamma\to\infty$.

When $\beta= 0$, some of the arguments fail, because $G(s, \gamma)$
might not be continuous at $(0,0)$. Assume that $F \in\operatorname
{AEP}(b, C)$
for some $b \geq2$ and $C > 0$. Note that $M = \infty$ by Lemma \ref
{lem:aep-rate}.
If $b = 2$, $G(s, \gamma) \to C \gamma$ when $\gamma> 0$ is fixed and
$s \to0$, by Lemma~\ref{lem:aep-rate}, and taking this limit as an extension
at $s = 0$, the same arguments used in the case $\beta> 0$ apply. If
$b > 2$, we need slightly different arguments.
As before, let $s_\gamma$ be a minimizer of $G(s,\gamma)$. We have that
$s_\gamma$ is well defined for all $\gamma$ and strictly positive,
because $G$ is uniformly continuous on any compact of $(0, 1/\zeta]
\times[0, \infty)$ and $G(s,\gamma) \sim C \gamma^{b/2} s^{1-b/2}
\to\infty$ when $s \to0$. Thus we may proceed as before in (\ref{L-ub})--(\ref
{L-lb}), obtaining that $L_0(\gamma)$ is strictly
increasing and continuous.
%
As before, we turn to proving that $L_0$ takes values $<$1 and finite
values $>$1.
First, with $\gamma= (\mu-\nu)^2/(2\zeta)$ and $s = 1/(2\zeta)$,
\[
L_0(\gamma) \leq G(s, \gamma) = \gamma\zeta/(\mu-\nu)^2 = 1/2 < 1.
\]
Next, showing that $L_0$ takes finite values above 1 is done exactly as
before, except that (\ref{v3}) is replaced by
\[
G(s, \gamma) \sim C s^{1-b/2} \gamma^{b/2} \geq C \zeta^{b/2-1}
\gamma
^{b/2},\qquad \gamma\to\infty
\]
by Lemma~\ref{lem:aep-rate}.
\end{pf}

The following result describes the variations of $\gamma$ (defined in
Lemma~\ref{lem:gamma}) with the parameter of an exponential family.
\begin{lem} \label{lem:gamma-prop}
Consider a natural exponential family of distributions $(F_\theta,
\theta\geq0)$ and let $\mu_\theta$ and $\Lambda_\theta^*$ denote the
mean and the rate function of $F_\theta$, respectively. Let $\zeta
_\theta$ be a continuous and decreasing function of $\theta$. Then, for
any fixed $0 < \beta< 1/\zeta_0$, $\gamma_\theta:= \gamma(F_\theta,
\mu_0, \zeta_\theta, \beta)$ is continuous and strictly increasing in
$\theta$. Moreover, if $\zeta_\theta\to0$ when $\theta\to\theta_c$,
then $\gamma_\theta\to\infty$ when $\theta\to\theta_c$.
\end{lem}

\begin{pf}
First, note that $\mu_\theta\geq\mu_0$ (Brown~\cite{MR882001}, Cor.
2.22) so
that $\gamma_\theta$ is well-defined. That $\gamma_\theta$ is strictly
increasing comes from the fact that both $\zeta_\theta$ and $\Lambda
_\theta^*(a)$ ($a > \mu_\theta$ fixed) are decreasing. The latter can
be seen from
\[
\Lambda_\theta^*(a) = - \lim_{k \to\infty} \frac1k \log\bbP
_\theta
(\bar{X}_k \geq a),
\]
where $\bar{X}_k$ is the average of the sample of size $k$ from
$F_\theta$ Brown~\cite{MR882001}, Cor. 2.22, and the fact that the
distribution of $\bar{X}_k$ as $\theta$ varies forms a natural
exponential family with parameter $k \theta$. That $\gamma_\theta$ is
continuous comes from the continuity of $\zeta_\theta$ and $\Lambda
_\theta^*(a)$ (in $(\theta, a)$).

For the behavior near $\theta_c$, note that $\Lambda_\theta^*(a) = 0$
for $a \leq\mu_\theta$, so that $G(1/(2\zeta_\theta), \gamma) = 1/2$
for any $\gamma
\leq(\mu_\theta-\mu_0)^2/(2\zeta_\theta)$. Combine this with the fact
that $\mu_\theta$ is strictly increasing in $\theta$ to see that
$\gamma
_\theta$ is of order at least $1/\zeta_\theta$. In fact, it is easy to
see that $\gamma_\theta\sim(\mu_\theta-\mu_0)^2/\zeta_\theta$ when
$\theta\nearrow\theta_c$.\vspace*{-3pt}
\end{pf}
%
\subsection{Main proofs}\vspace*{-3pt}
\label{sec:proofs}
%
\subsubsection{\texorpdfstring{Proof of Theorem \protect\ref{thm:cc-cube}}{Proof of Theorem 1}}
\label{sec:proof-cc-cube}
By monotonicity, it is sufficient to assume that $\theta_m = \theta$
for all $m$.
Fix $t$ and, for short, let $p = p_0(t)$ and $p' = p_\theta(t)$.
First, assume that $\theta> \theta_*$, so that $\zeta_{p'} < \alpha
\zeta_{p}$. Fix $B$ such that $1/\zeta_p < B < \alpha/ \zeta_{p'}$ and
consider the test with rejection region $\{S_m(t) \geq d B \log m\}$.
Under $\bbH^m_0$, we have $S_m(t) = (1+\o_\mathrm{P}(1)) (d/\zeta
_p) \log m$ by
(\ref{S-conv}), so that $\mathbb{P}(S_m(t) \geq d B \log m) \to0$.
Under $\bbH
^m_{1,K}$, $S_m(t) \geq S_K(t) = (1+\o_\mathrm{P}(1)) (\alpha d/\zeta
_{p'}) \log
m$, so that $\mathbb{P}(S_m(t) \geq d B \log m) \to1$. Thus this test is
asymptotically powerful.

Now assume that $\theta< \theta_*$, so that $\zeta_{p'} > \alpha
\zeta
_p$ and there is $B$ such that $\alpha/\zeta_{p'} < B < 1/ \zeta_p$.
Let $K^c = \Vm\setminus K$. It is sufficient to show that under both
$\bbH^m_0$ and $\bbH^m_{1,K}$, $S_m(t) = S_{K^c}(t)$ with probability
tending to 1, so that the values at the nodes in $K$ have no influence
on $S_m(t)$. Indeed, let $J$ be a hypercube within $\Vm$ of sidelength
$[m/3]$ which does not intersect $K$. Then $S_{K^c}(t) \geq S_J(t)$,
and the distribution of $S_{J}(t)$ is the same under both $\bbH^m_0$ and
$\bbH^m_{1,K}$. In addition, $\mathbb{P}(S_J(t) \geq d B \log m) \to
1$ by
(\ref{S-conv}).
Now, let $L$ be the set of nodes within (supnorm) distance $(\log m)^2$
from $K$, so that $L$ is a hypercube of side length $[m^\alpha] +
[2(\log m)^2]$ containing $K$ in its interior. In the event that $\{
S_m(t) \leq(\log m)^2\}$, $S_m(t) \neq S_{K^c}(t)$ only when $S_L(t) >
S_{K^c}(t)$. The distribution of $S_L(t)$ under the null is
stochastically bounded by its distribution under $\bbH^m_{1,K}$, which
is itself bounded by its distribution under $\bbH^m_{1,L}$. Even under
the latter, $\mathbb{P}(S_L(t) \geq d B \log m) \to0$ by (\ref
{S-conv}). We
then conclude the proof using the fact that $\mathbb{P}(S_m(t) \leq
(\log m)^2) \to1$, again by (\ref{S-conv}).\vspace*{-3pt}
%
\subsubsection{\texorpdfstring{Proof of Theorem \protect\ref{thm:cc-path}}{Proof of Theorem 2}}
\label{sec:proof-cc-path}
Here we use the notation and follow the arguments of Section \ref
{sec:proof-cc-cube}. In addition, let $\zeta^1_{p'} = \log(1/p')$, that is,
the function $\zeta$ in dimension one. When $\theta> \theta_*^+$, we
consider $1/\zeta_p < B < \alpha/d \zeta^1_{p'}$. Under $\bbH^m_0$, we
still have $S_m(t) = (1+\o_\mathrm{P}(1)) (d/\zeta_p) \log m$. Under
\vspace*{1pt}
$\bbH
^m_{1,K}$, $S_m(t) \geq S_K(t) = (1+\o_\mathrm{P}(1)) (\alpha/\zeta
_{p'}) \log
m$, because $K$ is isomorphic to a subinterval
of the one-dimensional lattice. We conclude as before that the test
with rejection region $\{S_m(t) \geq d B \log m\}$ is asymptotically powerful.

When $\theta< \theta_*^-$, we consider $\alpha/d\zeta_{p'} < B <
1/\zeta_p$. As before, let $L$ be the set of nodes within (supnorm)
distance $(\log m)^2$ from $K$, so that $L$ is now a band. As before,
it suffices to prove that $\mathbb{P}(S_L(t) \geq d B \log m) \to0$ under
$\bbH^m_{1,L}$. Although (\ref{S-conv}) cannot be applied, because $L$
is not isomorphic to a square lattice, its proof via the union bound
and (\ref{zeta}) applies. Indeed, fix $\eta> 0$ small enough that
$(1-\eta) \zeta_{p'} d B > \alpha$. Then, for $m$ large enough, we have
\begin{eqnarray*}
\mathbb{P}\bigl(S_L(t) \geq d B \log m\bigr)
&\leq&|L| \cdot\mathbb{P}(S \geq d B \log m) \\
&\leq&\O\bigl(m^{\alpha} (\log m)^{2(d-1)}\bigr) \exp\bigl(-(1-\eta) \zeta_{p'} d B
\log m\bigr) \\
&=& \O(\log m)^{2(d-1)} \exp\bigl( \bigl(\alpha- (1-\eta) \zeta_{p'} d B\bigr)
\log
m\bigr) \to0.
\end{eqnarray*}
%
\subsubsection{\texorpdfstring{Proof of Proposition \protect\ref{prp:t-large}}{Proof of Proposition 1}}
%
Let $k_m(\eps) = (1-\eps) d \log(m)/\log(1/p_0(t_m))$ with $\eps> 0$
fixed. We first show that $S_m(t_m) \geq k_m(\eps)$ with probability
tending to 1 under $\bbH^m_0$. We use the notation and arguments
provided in the proof of Lemma~\ref{lem:S-conv}. As in (\ref{Sm-lb}),
\begin{eqnarray*}
\mathbb{P}\bigl(S_m(t_m) < k_m(\eps)\bigr)
&\leq&\bigl(1 -\mathbb{P}\bigl(S \geq k_m(\eps)\bigr)\bigr)^N \\[-2pt]
&\leq&\bigl(1 -p_0(t_m)^{k_m(\eps)} \bigr)^N \\[-2pt]
&\leq&\exp\bigl(- m^{\eps d}/(\log m)^{2d} \bigr) \to0,
\end{eqnarray*}
where the second inequality holds for $m$ large enough by Lemma \ref
{lem:p-small}.

Assume that $\theta_m \leq\theta<\infty$ for all $m$. Proceeding as
in Section~\ref{sec:proof-cc-cube} and using the slightly larger
region $L$, it
is sufficient to show that for $\eps$ small enough, $S_L(t_m) \leq
k_m(\eps)$ when $X_v \sim F_\theta$ for all $v \in L$. Using the union
bound and the fact that $|L| = \O(m)^{\alpha d}$, we have\looseness=-1
%
\begin{equation}\label{SL-ub}
\mathbb{P}\bigl(S_L(t_m) \geq k_m(\eps)\bigr) \leq|L| \cdot\mathbb{P}\bigl(S \geq
k_m(\eps)\bigr) \leq
\O(m)^{\alpha d} (c p_\theta(t_m))^{k_m(\eps)},
\end{equation}
where the last inequality is due to Lemma~\ref{lem:p-small} (and $c$
is the
constant that appears there).
Through integration by parts, for $\theta> 0$ and $\eps\in(0,1)$ fixed,
we have $p_\theta(t) \leq p_0((1-\eps)t)$ for sufficiently large $t$.
Indeed, for $t$ large enough,
%
\begin{eqnarray*}
p_\theta(t)
&=& \exp\bigl(\theta t -\Lambda(\theta)\bigr) p_0(t) + \int_t^\infty\theta
\exp
\bigl(\theta x -\Lambda(\theta)\bigr) p_0(x)\,\mathrm{d}x \\[-2pt]
&\leq&\exp\bigl(\theta t -\Lambda(\theta) -C (1-\eps/3)^b t^b\bigr) + \int
_t^\infty\theta\exp\bigl(\theta x -\Lambda(\theta) -C (1-\eps/3)^b
x^b\bigr) \,\mathrm{d}x
\\[-2pt]
&\leq&\exp\bigl(-C (1-\eps/2)^b t^b\bigr) \\[-2pt]
&\leq& p_0\bigl((1-\eps) t\bigr),
\end{eqnarray*}
where we used the fact that $b > 1$ in line 3 and the fact that $\log
p_0(t) \sim-C t^b$ as $t \to\infty$ (because $F_0 \in\operatorname{AEP}(b,
C)$) in lines 2 and 4. The last property also implies that $p_0((1-\eps
) t) \leq p_0(t)^{(1-\eps)^{b+1}}$
for large $t$. Thus, for $m$ large enough,
$p_{\theta}(t_m) \leq p_0(t_m)^{(1-\eps)^{b+1}}$, so that taking logs
in (\ref{SL-ub}), we get
\[
\log\mathbb{P}\bigl(S_L(t_m) \geq k_m(\eps)\bigr) \leq\O(1) + (d \log m)
\bigl(\alpha
+ \O(\log p_0(t_m))^{-1} -(1-\eps)^{b+2} \bigr) \to-\infty,
\]
when $\eps< 1 -\alpha^{1/(b+2)}$. (Remember that $\alpha< 1$ and that
$p_0(t_m) \to0$, so the middle term is small.)
%
%
\subsubsection{\texorpdfstring{Proof of Theorem \protect\ref{thm:cc-cube-sup}}{Proof of Theorem 3}}
Let $\E_\theta$ denote the expectation of $X_v$ under $F_\theta$.
By Lemma~\ref{lem:large-sup}, under the null,
%
\begin{equation}\label{clt}
\frac{S_m(t) -\E_0(S_m(t))}{\sqrt{\Var_0(S_m(t))}} \to\cN(0,1),
\end{equation}
with $\Var_0(S_m(t))$ of order $m^d$.
Write $p:=p_0(t)$ and $p':= p_{\theta_m}(t)$.\vadjust{\goodbreak}

We consider the alternative with anomalous cluster $K$ as a
two-stage percolation process, where the first stage is percolation on
$\Vm$ with
probability $p$, as under the null, and the second stage is
percolation on the closed nodes within $K$, that is, $K
\setminus\{v\dvt
X_v > t\}$,
with (conditional) probability $(p' -p)/(1-p)$. An open cluster at the
first stage is called \emph{small} if it is not
a largest open cluster.

We may assume, except where noted below, that $\theta_m \to0$.
Because
\[
\frac{\partial}{\partial\theta} \log p_\theta(t) = \E_\theta(X_v
| X_v
> t) -\E_\theta(X_v),
\]
which is positive at $\theta= 0$ by choice of $t$, there exists $c\in
(0,\infty)$ such
that
%
\begin{equation}\label{theta2}
p'-p \sim c \theta_m \qquad\mbox{as } m \to\infty.
\end{equation}
%

Let $\Delta_m\ge0$ be the difference between the sizes of the largest
clusters under the null and the alternative.
For $x\in K$, let $F_x$ be the sum of the sizes of all small clusters
of the entire lattice
that contain some neighbor of $x$. Note that $\Delta_m \le\sum_{x\in
D}(1+F_x)$,
where $D$ is the set of $x \in K$ that are closed at the first stage
and open
at the second stage.
Therefore,
$\Delta_m$ has expectation bounded above by
%
\begin{equation}\label{delta3}
\E(\Delta_m) \le\biggl(\frac{p'-p}{1-p} \biggr) |K|(1+2d\mu_p),
\end{equation}
where $\mu_p < \infty$ is the mean size of a finite open cluster in the
infinite lattice.

By (\ref{theta2}) and the foregoing, $\E(\Delta_m) \le C \theta_m
m^{\alpha d}$ for some $C < \infty$.
By Markov's inequality,
$\Delta_m = \OP(\theta_m m^{\alpha d})$.

Thus, if $\theta_m m^{(\alpha-1/2) d} \to0$, then $\Delta_m/\sqrt
{\Var
_0(S_m(t))} \to0$, implying that the same central limit law as (\ref
{clt}) holds under the alternative, so
that the test based on the largest open cluster
is asymptotically powerless. We also must consider the case where
$\theta_m \not\to0$,
for which a similar argument is valid.

Now assume that $\alpha\geq1/2$ and $\theta_m m^{(\alpha-1/2) d}
\to
\infty$.
By Grimmett~\cite{MR1707339}, Theorem 8.99, and standard properties of the largest
cluster in a box (to be found in, e.g.,
Falconer and Grimmett~\cite{FG92}), with probability tending to 1, the
largest open cluster
increases in size by at least $C_1(p'-p) |K|$ for some $C_1=C_1(p)>0$.
By (\ref{theta2}), this has order $\theta_m m^{\alpha d}$.
Because
\[
\frac{\theta_m m^{\alpha d}}{\sqrt{\Var_0(S_m(t))}} \sim C_2 \theta_m
m^{(\alpha-1/2)d} \to\infty
\]
for some $C_2=C_2(p)>0$, the test based on the largest open cluster is
asymptotically powerful.
%
%
%
%
\subsubsection{\texorpdfstring{Proof of Theorem \protect\ref{thm:cc-cube-cri}}{Proof of Theorem 4}}
We may assume without loss of generality that $\theta_m \to0$ as $m
\to
\infty$.
By (\ref{S-cri}) and the assumption on $t_m$, we have that $S_m(t_m)
\asymP\log m$ under the null.
Now $p_\theta(t)$ is infinitely differentiable in $\theta$, with each
derivative continuous in $t$ and with
\[
\frac{\partial p_\theta(t)}{\partial\theta} \Big|_{\theta=
0} = p_0(t) [ \bbE_0(X_v|X_v > t) -\bbE_0(X_v) ] \geq\frac
{p_c}{2} [ \bbE_0(X_v|X_v > t_c) -\bbE_0(X_v) ] > 0,
\]
uniformly for $t$ in a neighborhood of $t_c$.
Therefore, there exists $C>0$ such that
\[
\frac{\partial p_\theta(t)}{\partial\theta} \geq1/C\quad
\mbox{and}\quad \biggl| \frac{\partial^2 p_\theta(t)}{\partial\theta^2}
\biggr| \leq C
\]
for $(\theta, t)$ in some neighborhood of $(0, t_c)$. Thus,
\[
p_\theta(t) -p_0(t) \geq\theta/C -C^2 \theta^2/2 \geq\theta/(2 C),
\]
on such a neighborhood. Let $A$ and $B$ be such that $p_c -p_0(t_m)
\leq A m^{-\alpha/\nu'}$ and $\theta_m \geq B m^{-\alpha/\nu'}$, and
assume that $B > 2 A C$, based on the statement of the theorem.
Because $\theta_m \to0$ and $t_m \to t_c$,
\[
m^{\alpha/\nu''} \bigl(p_{\theta_m}(t_m) -p_c\bigr) \geq
m^{\alpha/\nu''} \biggl[\frac{\theta_m}{2C} + \bigl(p_0(t_m)-p_c\bigr) \biggr] \geq
\biggl[\frac{B}{2C} -A \biggr] m^{\alpha(1/\nu'' -1/\nu')} \to\infty
\]
for $\nu''<\nu'$ and sufficiently large $m$.
By (\ref{S-cri}) applied to $K \in\Km$, it follows that $S_K(t_m)
\asymP m^{\alpha d}$ under the alternative.
Consequently, the test with rejection region $\{S_m(t_m) \geq(\log
m)^2\}$
is asymptotically powerful.
%
\subsubsection{\texorpdfstring{Proof of Lemma \protect\ref{lem:uls-sub}}{Proof of Lemma 5}}
\textit{Part} 1. This follows immediately from Lemma~\ref{lem:S-conv}.

Therefore, we focus on the remaining two parts.
We use the abbreviated notation $F := F_{\theta|t}$, $\Lambda^* :=
\Lambda^*_{\theta|t}$, $\mu:= \mu_{\theta|t}$, $\zeta:= \zeta
_{p_\theta(t)}$, $\gamma:= \gamma_{\theta|t}(\beta)$, $U_m := U
_m(t, k_m)$, and write $\nu:= \mu_{0|t}$.
Let $Y_k=X_k-\nu$.
%
%
As in Lemma~\ref{lem:Nk}, let $N_m(k)$ denote the number of open
cluster of
size $k$ within $\Vm$, and define
\[
G_k(x) = \mathbb{P}(k^{1/2} \bar{Y}_{k} \leq x),
\]
where $\bar{Y}_{k} = \bar{X}_{k} -\nu$ and $\bar{X}_k$ is the average
of an i.i.d. sample of size $k$ from $F$.
By the independence of $\bar{Y}_K$ and $\bar{Y}_L$ for $K, L \in\Qmt$
distinct, we have
\[
\mathbb{P}(U_m \leq x) = \mathbb{E}\biggl(\prod_{k \geq k_m}
G_k(x)^{N_m(k)}\biggr) =
\mathbb{E}(\exp[-R_m(x)]),
\]
where
\[
R_m(x) := - \sum_{k \geq k_m} N_m(k)\log\bigl(1 -\bar{G}_k(x)\bigr).
\]
Thus, we turn to bounding $R_m(x)$.

\textit{Part} 2. Define $x_m = \sqrt{\gamma d \log m}$ and fix $\eps
> 0$.
For the lower bound, let $\ell_m$ be the closest integer to $a d \log
m$ between $k_m$ and $(d/\zeta) \log m$, where
%
\begin{equation}\label{a}
a = \mathop{\argmin}_{\beta< s < 1/\zeta} \bigl[s \Lambda^* \bigl(\nu+\sqrt
{\gamma/s} \bigr) + s \zeta\bigr].
\end{equation}
We have
\[
R_m\bigl((1-\eps) x_m\bigr)
\geq T_m := N_m(\ell_m) \bar{G}_{\ell_m}\bigl((1-\eps) x_m\bigr),
\]
and we show that for $\eps$ fixed, $T_m \to\infty$ in probability. Fix
$\eta> 0$.
On the one hand, we use Lemma~\ref{lem:Nk} and (\ref{zeta}), to get
\[
\mathbb{E}(N_m(\ell_m)) \geq\frac{(m -2\ell_m)^d}{\ell_m} \mathbb
{P}(S = \ell_m)
\geq m^d \exp\bigl(-(1+\eta) \zeta\ell_m\bigr)
\]
for $m$ large enough.
On the other hand, we use Cram\'er's theorem (Dembo and Zeitouni~\cite
{MR2571413}, Theorem 2.2.3) to get
\begin{eqnarray*}
\bar{G}_{\ell_m}\bigl((1-\eps) x_m\bigr)
&\geq&\mathbb{P}\bigl(\bar{Y}_{\ell_m} \geq(1-\eps/2) \sqrt{\gamma
/a}\bigr) \\
&\geq&\exp\bigl(-(1+\eta) \ell_m \Lambda^* \bigl[\nu+(1-\eps/2) \sqrt
{\gamma/a} \bigr]\bigr )
\end{eqnarray*}
for $m$ large enough.
By the definition of $\gamma$, $a \Lambda^* [\nu+\sqrt{\gamma
/a} ] + a \zeta= 1$, and thus for $\eps$ small enough,
\[
a \zeta+ a \Lambda^* \bigl[\nu+(1-\eps/2) \sqrt{\gamma/a} \bigr] < 1,
\]
by strict monotonicity, as in the proof of Lemma~\ref{lem:gamma}.
Thus, for $\eta$ small enough,
\[
\ell_m \zeta+ \ell_m \Lambda^* \bigl[\nu+(1-\eps/2) \sqrt{\gamma
/a} \bigr] \leq(1-\eta) d \log m.
\]
It follows that
\[
\mathbb{E}(T_m) \geq m^{\eta^2 d}.
\]
To bound the corresponding variance, we use Lemma~\ref{lem:Nk} to obtain
\[
\operatorname{Var}(T_m) \leq\O(\log m)^{d} \mathbb{E}(T_m),
\]
and it follows by Chebyshev's inequality that indeed $T_m \to\infty$
in probability.

Because $T_m\ge0$, $\exp(-T_m) \to0$ in $L^1$, and thus
\[
\mathbb{P}\bigl(U_m\le(1-\eps)x_m\bigr)\to0.
\]

We next show that $\mathbb{E}(R_m((1+\eps) x_m)) \to0$, which will imply
the claim of Part 2.
Fix $\eta> 0$.
We have that
%
\begin{equation}\label{oh1}
R_m\bigl((1+\eps) x_m\bigr) \leq T_m +2Z_m,
\end{equation}
where
\[
T_m := 2 \sum_{k = k_m}^{k_m^{(\eta)}} N_m(k) \bar{G}_k\bigl((1+\eps) x_m\bigr)
\]
and $Z_m$ is the number of clusters of size exceeding
$k_m^{(\eta)} := [(1+\eta) (d/\zeta) \log m]$.
We first note that, as in the proof of Lemma~\ref{lem:Nk}, for large $m$,
%
\begin{equation}\label{oh2}
\mathbb{E}(Z_m) \le m^d \exp\bigl(-\tfrac12\zeta k_m^{(\eta)}\bigr)\to0.
\end{equation}

We next turn to $T_m$, and show that for $\eps$ fixed and $\eta$ small
enough, $\mathbb{E}(T_m) \to0$.
On the one hand, we use Lemma~\ref{lem:Nk} and (\ref{zeta}) to get
\[
\mathbb{E}(N_m(k)) \leq m^d \exp\bigl(-(1-\eta) \zeta k\bigr)
\]
for $m$ large enough.
On the other hand, by Chernoff's bound,
\[
\bar{G}_{k}\bigl((1+\eps) x_m\bigr)
\leq\exp\bigl(-k \Lambda^*\bigl [\nu+(1+\eps) x_m/\sqrt{k} \bigr]\bigr ).
\]
Taken together, we obtain
\begin{eqnarray*}
\mathbb{E}(T_m)
&\leq&2 \sum_{k = k_m}^{k_m^{(\eta)}} m^d \exp\bigl(-(1-\eta)\bigl [k
\zeta+ k \Lambda^* \bigl(\nu+(1+\eps) x_m/\sqrt{k}\bigr )\bigr ] \bigr)
\\
&\leq&\O(\log m) \exp\Bigl(d \log m -(1-\eta)\min_{k_m \leq k \leq
k_m^{(\eta)}} \bigl[k \zeta+ k \Lambda^* \bigl(\nu+(1+\eps) x_m/\sqrt
{k} \bigr) \bigr] \Bigr) \\
&\leq&\O(\log m) \exp\bigl( \bigl(1 - (1-\eta) A\bigr) d \log m\bigr ),
\end{eqnarray*}
where
%
\begin{equation}\label{A}
A := \inf_{\beta< a < (1+\eta)/\zeta} \bigl[a \Lambda^* \bigl(\nu
+(1+\eps)\sqrt{\gamma/a}\bigr ) + a \zeta\bigr].
\end{equation}
As in the proof of Lemma~\ref{lem:gamma}, $A = A(\eps, \eta)$ is
continuous in
$(\eps, \eta)$ and strictly increasing in $\eps$. Because $A(0,0) = 1$
by definition of $\gamma$, for $\eps$ fixed, $-h := 1 - (1-\eta)
A(\eps
, \eta) < 0$ for $\eta$ small enough, in which case $\mathbb{E}(T_m)
\leq
m^{-hd/2} \to0$ as $m$ increases.

By (\ref{oh1})--(\ref{oh2}), we have that $\mathbb{E}(R_m((1+\eps
)x_m)))\to0$. By Jensen's inequality,
\[
\mathbb{P}\bigl(U_m \le(1+\eps)x_m\bigr) \ge\exp\bigl(-\mathbb{E}\bigl(R_m\bigl((1+\eps
)x_m\bigr)\bigr)\bigr) \to1,
\]
and the proof of this part is complete.
%

\textit{Part} 3. We build on the arguments provided so far, which apply
essentially unchanged, except in two places. In the lower bound,
instead of Cram\'er's theorem, we use
\[
\bar{G}_k(x) \geq\bar{F}\bigl(x/\sqrt{k}\bigr)^k,
\]
combined with the asymptotic behavior for $\bar{F}$. In the upper
bound, $A$ defined in (\ref{A}) is evaluated differently when $b < 2$.\vadjust{\goodbreak}

\textit{Part} 3(a). When $b > 2$, we have $a > 0$ in (\ref{a}) (with
$\beta= 0$), because
\[
h(s) := s \Lambda^* \bigl(\nu+\sqrt{\gamma/s} \bigr) + s \zeta\asymp
s^{1-b/2} \to\infty
\]
for $\gamma$ fixed and $s \to0$, by Lemma~\ref{lem:aep-rate}. When
$b = 2$,
we take $a$ small enough if the minimum is at $a = 0$. Then the other
arguments in Part 2 apply unchanged.

\textit{Part} 3(b). By the same calculations, $a = 0$ in (\ref{a}),
because $h(s) > 0$ for all $s > 0$, and $h(s) \asymp s^{1-b/2} \to0$
when $s \to0$, because $b < 2$. This would make $A = 0$ in (\ref{A})
for any $\eps> 0$, making the arguments for the upper bound collapse.
Instead, redefine $x_m = (C d \log m)^{1/b} k_m^{1/2 -1/b}$. Because
$x_m/\sqrt{k} \to\infty$ uniformly over $k \leq k_m^{(\eta)}$, for
$\eta> 0$ fixed, we have
\[
k \zeta+ k \Lambda^* \bigl(\nu+(1+\eps) x_m/\sqrt{k}\bigr ) \geq k
\zeta+ (1-\eta) C k^{1-b/2} (1+\eps)^b x_m^b
\]
for $m$ large enough, by Lemma~\ref{lem:aep-rate}. Then the term on the
right-hand side takes its minimum over $k_m \leq k \leq k_m^{(\eta)}$
at $k = k_m$, and from here, the remaining arguments apply.
%
%
\subsubsection{\texorpdfstring{Proof of Proposition \protect\ref{prp:uls-bad}}{Proof of Proposition 2}}
Assume, for simplicity, that $\theta_m = \theta< \theta_c$ for all $m$.
The key point is that $F_{\theta|t} \in\operatorname{AEP}(b, C)$.
Indeed, we
have $\bar{F}_{\theta|t}(x) = \bar{F}_{\theta}(x)/\bar{F}_{\theta}(t)$,
where the denominator is constant in $x$ and, integrating by parts,
\[
\bar{F}_{\theta}(x) = \exp\bigl(\theta x -\Lambda(\theta)\bigr) \bar
{F}_{0}(x) +
\int_x^\infty\theta\exp\bigl(\theta y -\Lambda(\theta)\bigr) \bar
{F}_{0}(y)\,\mathrm{d}y.
\]
From here, we reason as in the proof of Proposition~\ref{prp:t-large},
using the
fact that $\log\bar{F}_{0}(y) \sim-C y^b$ when $y \to\infty$, with
$b > 1$.
Thus $F_{\theta|t}$ and $F_{0|t}$ have same (first-order) asymptotics,
and so nothing distinguishes the asymptotic behavior of $U_m$ under
the null and under an alternative.
In detail, we proceed as in Section~\ref{sec:proof-cc-cube}, with the enlarged
hypercube $L$, and show that in probability under $\bbH^m_{1,L}$,
\[
\limsup_{m \to\infty} k_m^{1/b -1/2} (\log m)^{-1/b} U_L < (d/C)^{1/b},
\]
where $U_L$ is the $\operatorname{ULS}$  scan statistic restricted to open clusters
within $L$. Because $L$ is a scaled version of $\Vm$, $F_{\theta|t}
\in
\operatorname{AEP}(b, C)$ and $p_\theta(t) < p_c$, Lemma~\ref{lem:uls-sub}
applies to yield
\[
k_m^{1/b -1/2} (\alpha\log m)^{-1/b} U_L \to(d/C)^{1/b}.
\]
We then conclude with the fact that $\alpha< 1$.
%
\subsubsection{\texorpdfstring{Proof of Theorem \protect\ref{thm:uls-cube} and Theorem \protect\ref{thm:uls-path}}
{Proof of Theorem 5 and Theorem 6}}
The proof of Theorem~\ref{thm:uls-cube} is parallel to that of Theorem
\ref{thm:cc-cube}
in Section~\ref{sec:proof-cc-cube}, but using Lemma~\ref{lem:uls-sub}
in place of
Lemma~\ref{lem:S-conv}.
Note that we use the fact that for $t$ and $\beta> 0$ fixed, $\gamma
_{\theta|t}(\beta)$ is continuous and strictly increasing in $\theta$.
This comes from Lemma~\ref{lem:gamma-prop} and the fact that when $t$ is
fixed, $F_{\theta|t}$ is also a natural exponential family with
parameter~$\theta$.
Similarly, the proof of Theorem~\ref{thm:uls-path} is parallel to that
of Theorem~\ref{thm:cc-path} in Section~\ref{sec:proof-cc-path}.
Further details are omitted.
%
\subsubsection{\texorpdfstring{Proof of Lemma \protect\ref{lem:uls-sup}}{Proof of Lemma 6}}
The proof is parallel to that of Lemma~\ref{lem:uls-sub}. In
particular, we
use the notation introduced there and only note where the arguments
differ (although never substantially).

\textit{Part} 1. In this case, by Lemma~\ref{lem:large-sup} and Lemma
\ref{lem:small-sup}, there is only one open cluster with size $k_m$ or larger,
and the result follows from, for example, Chebyshev's inequality.

\textit{Part} 2. Define $x_m = \sqrt{2 \sigma^2 d(1 -\delta\beta')
\log
m}$ and fix $\eps> 0$. For the lower bound, we have
\[
R_m\bigl((1-\eps) x_m\bigr)
\geq T_m := N_m(k_m) \bar{G}_{k_m}\bigl((1-\eps) x_m\bigr).
\]
Fix $\eta> 0$.
By Lemma~\ref{lem:Nk} (still valid) and (\ref{delta}),
\[
\mathbb{E}(N_m(k_m)) \geq m^d \exp\bigl(-(1+\eta) \delta k_m^{(d-1)/d}\bigr)
\]
for $m$ large enough.
By Cram\'er's theorem and the fact that $\Lambda^*(x) \sim x^2/(2
\sigma
^2)$ when $x$ is small,
\begin{eqnarray*}
\bar{G}_{k_m}\bigl((1-\eps) x_m\bigr)
&\geq&\exp\bigl(-(1+\eta) k_m \Lambda^* \bigl[(1-\eps) x_m/\sqrt
{k_m} \bigr]\bigr ) \\
&\geq&\exp\bigl(-(1+\eta) (1-\eps/2) x_m^2/(2 \sigma^2) \bigr)
\end{eqnarray*}
for $m$ large enough.
Thus,
\[
\mathbb{E}(T_m) \geq\exp\bigl(d \log m -(1+\eta) \bigl(\delta k_m^{(d-1)/d} +
(1-\eps/2) x_m^2/(2 \sigma^2)\bigr)\bigr ) \geq m^{\eps d (1 - \delta\beta')/4}
\]
for $m$ large enough and $\eta$ small enough.
For the variance, we use Lemma~\ref{lem:Nk} to get
\[
\operatorname{Var}(T_m) \leq\O(\log m)^{d^2/(d-1)} \mathbb{E}(T_m).
\]
We then conclude by Chebyshev's inequality.

We now show that $R_m((1+\eps) x_m) \to0$ in probability. Equation
(\ref{oh1}) holds
with $k_m^{(\eta)} := [(1+\eta) (d/\delta) \log m]^{d/(d-1)}$. As before,
\[
\mathbb{E}(Z_m) \le m^d \exp\bigl\{-\tfrac12\delta\bigl(k_m^{(\eta
)}\bigr)^{(d-1)/d}\bigr\} \to0 \qquad\mbox{as } m\to\infty.
\]
By Lemma~\ref{lem:Nk} and (\ref{delta}),
\[
\mathbb{E}(N_m(k)) \leq m^d \exp\bigl(-(1-\eta) \delta k^{(d-1)/d}\bigr)
\]
for $m$ large enough.
The absence of a boundary to $\Vm$ is
being used here. The tail behavior of percolation clusters near the
boundary of
a box is not yet fully understood (see the remark in Section \ref
{sec:uls-sup}).
By Chernoff's bound and the behavior of $\Lambda^*$ near the origin,
\[
\bar{G}_{k}\bigl((1+\eps) x_m\bigr)
\leq\exp\bigl(-(1+\eps) x_m^2/(2\sigma^2) \bigr)
\]
for any $k \geq k_m$.
Thus,
\begin{eqnarray*}
\mathbb{E}(T_m)
&\leq&2 \sum_{k = k_m}^{k_m^{(\eta)}} m^d \exp\bigl(-(1-\eta) \delta
k^{(d-1)/d} - (1+\eps) x_m^2/(2\sigma^2) \bigr) \\
&\leq&\O(\log m)^{d/(d-1)} m^{- \eps d (1 - \delta\beta')/4}
\end{eqnarray*}
for $m$ large enough and $\eta$ small enough.

\textit{Part} 3. This part is even more similar to what we did in the
proof of Lemma~\ref{lem:uls-sub}. The behavior of $U_m$ is driven by
the open
clusters of size of order $\log m$, with the only difference being that
the term in $k^{(d-1)/d}$ from the bounds on $N_m(k)$ is negligible.
Details are omitted.
%
\subsubsection{\texorpdfstring{Proof of Theorem \protect\ref{thm:uls-sup}}{Proof of Theorem 7}}
%
Without loss of generality, we assume that $\theta_m$ is bounded.
By Lemma~\ref{lem:uls-sup} and our assumptions on $k_m$, under the
null, $U_m
:= U_m(t, k_m) \simP A (\log m)^{1/2}$ for a finite constant $A > 0$.
We now consider the alternative, where the anomalous cluster is $K$.

The contribution of the largest open cluster, $Q_m$, is
\begin{eqnarray*}
\sqrt{|Q_m|} (\bar{X}_{Q_m} -\mu_{0|t})
&=& \frac{|Q_m \cap K|}{\sqrt{|Q_m|}} (\bar{X}_{Q_m \cap K} -\mu
_{\theta
_m|t}) + \frac{|Q_m \cap K^c|}{\sqrt{|Q_m|}} (\bar{X}_{Q_m \cap K^c}
-\mu_{0|t}) \\
&&{} + \frac{|Q_m \cap K|}{\sqrt{|Q_m|}} (\mu_{\theta_m|t}
-\mu_{0|t}).
\end{eqnarray*}
On the right-hand side, the first term is of order $\o_\mathrm
{P}(1)$, and the
second term is of order $\OP(1)$, by Chebyshev's inequality and the
fact that, with probability tending to 1, $|Q_m \cap K| \asymp|K|$
and $|Q_m| \asymp|\Vm|$, by Lemma~\ref{lem:large-sup}.
The last term is of (exact) order $\O(\theta_m m^{(\alpha-1/2)d})$,
by the fact that $\mu_{\theta|t}$ is differentiable at $\theta= 0$
with derivative equal to $\sigma_{0|t}^2 > 0$. Therefore, the $\operatorname{ULS}$
scan test is asymptotically powerful when
$\liminf\theta_m m^{(\alpha-1/2)d} (\log m)^{-1/2}$ is large enough.
(Note that this requires $\alpha> 1/2$.)
If instead, we have $\limsup\theta_m m^{(\alpha-1/2)d} (\log
m)^{-1/2} \to0$,
then the scan over $Q_m$ may be ignored, and we need to consider
smaller clusters.

By Lemma~\ref{lem:small-sup} and the upper bound on $k_m$, the second-largest
cluster entirely within $K$ is scanned and its contribution is of order
$\O(\theta_m (\log m)^{d/(2d-2)})$, by the same arguments that
established the contribution of the largest open cluster.
Thus, the $\operatorname{ULS}$  scan test is asymptotically
powerful when $\liminf\theta_m (\log m)^{d/(2d-2) -1/2}$ is large
enough.
If instead, $\theta_m (\log m)^{d/(2d-2) -1/2} \to0$,
the test is asymptotically powerless. Indeed, let $L$ be the set of
nodes within
distance $(\log m)^3$ from $K$, and let $U_L$ be the result of scanning
the open clusters of size at least $k_m$ and entirely within $L$.
As argued in the proof of Proposition~\ref{prp:uls-bad}, this time using
Lemma~\ref{lem:small-sup}, it is sufficient to show that
$U_L \leq A (\log m)^{1/2}$ with probability tending to 1 under
$\bbH_{1,L}^m$.
For any open cluster $Q$ entirely within $L$,
\[
\sqrt{|Q|} (\bar{X}_{Q} -\mu_{0|t})
= \sqrt{|Q|} (\bar{X}_{Q} -\mu_{\theta_m|t}) + \sqrt{|Q|} (\mu
_{\theta
_m|t} -\mu_{0|t}),
\]
so that
\[
U_L \leq\max_Q \sqrt{|Q|} (\bar{X}_{Q} -\mu_{\theta_m|t}) + \o
_\mathrm{P}(1),
\]
where the maximum is over open clusters of size at least $k_m$ and
entirely within $L$, and the second term is $\o_\mathrm{P}(1)$ by
Lemma~\ref{lem:small-sup} and the size of $\theta_m$.
Although $\theta_m \to0$ varies, this maximum may be handled exactly
as in Lemma~\ref{lem:uls-sup}, so that it is $\simP A (\alpha\log m)^{1/2}$,
and we conclude.
%
\subsubsection{\texorpdfstring{Proof of Lemma \protect\ref{lem:uls-orig}}{Proof of Lemma 7}}
%
%
%
We prove only the more refined part. We use abbreviated notation as
before, in particular, we omit the subscript 0, using $F_t = F_{0|t}$,
$\sigma_t = \sigma_{0|t}$, and so on. The lower bound is obtained via
$\ULS_m \geq U_m(t^*)/\sigma_{t^*}$, where $t^*$ defines $\Gamma
(\beta
)$, and applying Lemmas~\ref{lem:uls-sub} or~\ref{lem:uls-sup} to
$U_m(t^*)$ depending on whether $t^* > t_c$ or $t^* < t_c$. For
simplicity, we assume that $t^* \neq t_c$. If $t^* = t_c$, then we
consider a nearby threshold and argue by continuity. For the upper
bound, we prove that $\mathbb{P}(\ULS_m \ge x_m) \to0$, where $x_m
:= \sqrt{g
\log m}$ and $g > G:= (d \Gamma(\beta))^{1/2}$.

As $t$ increases, clusters are created and then destroyed in the
coupled percolation processes. Suppose the removal at time $t$ from the
percolation process of vertex $v$ creates some cluster $Q_t(w)$ at some
neighbor $w$ of $v$. If $\ULS_m \ge x_m$, there must exist a vertex $v$
and a neighbor $w$ such that the cluster formed at $w$ at time $X_v$
contributes at some future time $t' > X_v$ an amount at least $x_m$ to
$\ULS_m$. By conditioning on $v$, $X_v$, and $w$, one obtains that
%
\begin{equation}\label{final4}
\mathbb{P}(\ULS_m \ge x_m) \le\o(1) + \int_{-\infty}^{t_\beta}
\mathbb{P}\biggl(\bigcup_{v\in\Vm}\bigcup_{w \in\partial v} \Omega
_t(w)\biggr) \,\mathrm{d}F(t),
\end{equation}
where the $\o(1)$ term covers the probability that the cluster at time
$-\infty$, namely $\bbV_m$, determines $\ULS_m$, or that a cluster at
threshold $t > t_\beta$ is of size at least $k_m := \beta\log m$;
$\partial v$ is the neighbor set of $v$; and $\Omega_t(w)$ is the
event that:
\begin{enumerate}
\item$k := |Q_t(w)|$ satisfies $k \ge\beta\log m$,
\item there exists a time $t' \ge t$ such that $Q_t(w)$ still exists at
time $t'$ and
\item$Y_t(k)-\mathbb{E}(Y_{t'}(k)) \ge x_m \sigma_{t'}\sqrt k$, where
$Y_t(k)$ is the sum of a $k$-sample from $F_t$.
\end{enumerate}

Assume (briefly) that $\sigma_t$ is non-decreasing, and note that $\mu
_t$ is automatically non-decreasing. Then as in the proofs of
Lemmas~\ref{lem:uls-sub} and~\ref{lem:uls-sup}, and using similar notation,
\begin{eqnarray*}
\sum_{v\in\Vm}\sum_{w \in\partial v} \mathbb{P}(\Omega_t(w))
&\le&\sum_{v\in\Vm}\sum_{w \in\partial v} \mathbb{P}\bigl(k :=|Q_t(w)|
\ge\beta \log m, Y_t(k) - \mathbb{E}(Y_t(k)) \ge x\sigma_t\sqrt k\bigr)
\\
&\le&2d\E(R_{t}(x_m)), \qquad R_{t}(x) := \sum_{k \geq k_m} N_t(k)
\bar{G}_{t}(k,x),
\end{eqnarray*}
where $N_t(k)$ is the number of $t$-open clusters of size $k$ and
\[
\bar{G}_{t}(k,x) = \mathbb{P}\bigl(Y_t(k) - \mathbb{E}(Y_t(k)) \ge
x\sigma_t\sqrt k\bigr).
\]
Therefore, by (\ref{final4}),
%
\begin{eqnarray}\label{ULS-ub1}
 \mathbb{P}(\ULS_m \ge x_m) &\le&\o(1) + 2d \biggl( \int_{-\infty}^{t_c-h}
+ \int
_{t_c+h}^{t_\beta} \E(R_{t}(x_m))\,\mathrm{d}F(t) \biggr)\nonumber
\\[-8pt]
\\[-8pt]
&&{}+ F(t_c+h) - F(t_c-h)\nonumber
\end{eqnarray}
for any $h > 0$.
We bound $\E(R_{t}(x_m))$ as we did in the proofs of Lemmas \ref
{lem:uls-sub} and~\ref{lem:uls-sup}. Explicitly, when $t_c + h \leq t
\leq t_\beta$, we use Lemma~\ref{lem:Nk} and (\ref{S-ineq}), to get
\begin{eqnarray*}
\E(N_t(k))
&\leq&\bigl(1-p(t)\bigr)^2 \frac{k \mathrm{e}^{-k\zeta_{p(t)}}}{
(1 - \mathrm{e}^{-\zeta_{p(t)}})^2}
\\
&\leq&C(h,\beta) k \exp\bigl(-k \zeta_{p(t_c+h)}\bigr), \qquad
C(h,\beta) :=
\frac{(1-p(t_\beta))^2 }{(1 - \mathrm{e}^{-\zeta_{p(t_c+h)}})^2}.
\end{eqnarray*}
We use Chernoff's Bound on $\bar{G}_{t}(k,x)$, to obtain
\[
\E(R_{t}(x_m)) \leq C(h,\beta) (k_{m,t}^h)^2 \exp\bigl((1 - A_t) d \log
m\bigr ) + \exp\bigl(-h d \log(m)/2\bigr),
\]
where $k_{m,t}^h := (1+h)(d/\zeta_{p(t)}) \log m$,
\[
A_t := \inf_{\beta< s < (1+h)/\zeta_{p(t)}} \bigl[s \Lambda_t^*
\bigl(\mu+\sqrt{g/s} \bigr) + s \zeta_{p(t)} \bigr],
\]
as in (\ref{A}), and the last term is the probability that a there is a
$t$-open of size exceeding $k_{m,t}^h$.
Note that $A_t > 1$ for all $t_c + h \leq t \leq t_\beta$ because $g >
G$. By continuity of $A_t$, $A_+ := \inf\{A_t \dvt t_c + h \leq t \leq
t_\beta\} > 0$. Hence, we have the following bound for all $t_c + h
\leq t \leq t_\beta$,
\[
\E(R_{t}(x_m)) \leq C(h,\beta) \bigl[(1+h)\bigl(d/\zeta_{p(t_c+h)}\bigr) \log m\bigr]^2
m^{-(A_+ -1) d} + \exp\bigl(-h d \log(m)/2\bigr).
\]
When $t \leq t_c-h$, we simply use the fact that
\[
\sum_k \E(N_t(k)) \leq|\Vm| = m^d,
\]
and bound $\bar{G}_{t}(k,x)$ in the same way.
We get
\[
\E(R_{t}(x_m)) \leq\exp\bigl((1 - A_t) d \log m \bigr),
\]
where
\[
A_t := \inf_{\beta< s} s \Lambda_t^* \bigl(\mu+\sqrt{g/s} \bigr).
\]
Again, $A_t > 1$ for $t < t_c-h$ and $A_t \to A_{-\infty} > 1$ as $t
\to-\infty$. Hence, by continuity of $A_t$, $A_- := \inf\{A_t : t <
t_c -h\} > 0$, so that
\[
\E(R_{t}(x_m)) \leq m^{-(A_- -1) d},
\]
valid for all $t < t_c-h$.
Hence, the two integrals in (\ref{ULS-ub1}) tend to zero with $m$. We
then let $h \to0$ so that $F(t_c+h) - F(t_c-h) \to0$, because $F$ is
continuous at $t_c$.

Assume now that $F$ has no atoms on $(-\infty, t_\beta]$. Then
$\sigma
_{t}$ is continuous on $(-\infty, t_\beta]$, and in fact, is uniformly
continuous because $\sigma_{t} \to\sigma$ when $t \to-\infty$,
because it is positive on that interval (because $\sigma_t = 0$ implies
that $F_t$ is a point mass), $\ul\sigma:= \min\{\sigma_{t} \dvt t \leq
t_\beta\} > 0$. Because $g>G$ we can find $c>0$ such that
$g' := g(1-c)^2 >G$, and also $\eta>0$ such that
%
\begin{equation}\label{final6}
|\sigma_{s}-\sigma_{t}| \le c\ul\s, \qquad\mbox{if } |s-t|\le\eta,
s,t \leq
t_\beta.
\end{equation}
%
Let $x_m'=\sqrt{g'\log m}$.
We say that a cluster $Q$ \textit{scores at time} $s$ if it exists at time
$s$ and
in addition
\[
|Q| \ge\beta\log m, \qquad\sum_{v \in Q} X_v \ge|Q| \mu_s + x_m \sigma
_s \sqrt{|Q|}.
\]

Without loss of generality, assume that $t_c$ is not an integer
multiple of $\eta$.
Fix two neighbors $v,w \in\Vm$, and a time $t \leq t_\beta$. If
$\Omega
_t(w)$ occurs then either:
\begin{enumerate}[(b)]
\item[(a)] $Q_t(w)$ scores at some time $s \in[t, n_t \eta]$, where $n_t
\in\bbZ$ satisfies $(n_t-1)\eta\le t< n_t\eta$, or
\item[(b)] there exists $n \ge n_t$ and $s \in[n\eta, (n+1)\eta)$
such that $Q_{n\eta}(w)$ scores at time $s$.
\end{enumerate}
The latter possibility arises when $Q_t(w)$ scores at some
time $s$ not belonging to the interval $[t,n_t\eta)$. Writing
$[n\eta, (n+1)\eta)$ for the interval containing $s$, $Q_t(w)$ must
exist at the start of this interval, which is to say that
$Q_t(w)=Q_{n\eta}(w)$.

The probability of (a) is no larger than
%
\begin{equation}\label{final8}
\mathbb{P}\bigl(k := |Q_t(w)| \ge\beta\log m, \exists s\in[t,n_t\eta] \dvt
Y_t(k)/k \ge\mu_s + x_m \sigma_{s}/\sqrt k\bigr).
\end{equation}
By (\ref{final6}) and the fact that $\mu_s$ is non-decreasing,
%
\begin{equation}\label{g+3}
\mu_s+\frac{x_m \sigma_{s}}{\sqrt k} \ge\mu_t+ \frac{x_m'\sigma
_{t}}{\sqrt k},
\end{equation}
so that (\ref{final8}) is no greater than
%
\begin{equation}\label{case-a}
\mathbb{P}\bigl(k := |Q_t(w)| \ge\beta\log m, Y_t(k)/k \ge\mu_t + x_m'
\sigma _{t}/\sqrt k\bigr).
\end{equation}
%

Arguing similarly, part (b) has probability no greater than
%
\begin{equation}\label{case-b}
\sum_{t/\eta< n < t_\beta/\eta} \mathbb{P}\bigl(k := |Q_t(w)| \ge\beta
\log m, Y_{n\eta}(k)/k \ge\mu_{n\eta} + x_m'\sigma_{n\eta}/\sqrt k\bigr).
\end{equation}

We divide the integral in (\ref{final4}) as follows
\[
\int_{-\infty}^{t_\beta} = \int_{-\infty}^{-1/h} + \int
_{-1/h}^{t_c -h}
+ \int_{t_c-h}^{t_c+h} + \int_{t_c+h}^{t_\beta}.
\]
The first integral is bounded by $F(-1/h)$ and the third integral by
$F(t_c+h) - F(t_c-h)$, both terms vanishing as $h \to0$. For the
second and fourth integrals, we do exactly as before, separately for
(\ref{case-a}) and (\ref{case-b}) -- for the latter, the sum has at most
$(t_\beta+ 1/h)/\eta+ 1$ terms in the second integral and at most
$(t_\beta- t_c - h)/\eta+ 1$ terms in the fourth integral.
%
\subsubsection{\texorpdfstring{Proof of Theorem \protect\ref{thm:uls-orig}}{Proof of Theorem 8}}
By Lemma~\ref{lem:uls-orig}, $\ULS_m(k_m)$ is of order at most $\sqrt
{\log m}$
under the null.
Now consider the alternative with anomalous cluster $K$. If $0 <
(\alpha-1/2)d < \alpha/\nu$, consider the contribution of the largest
open cluster at supercritical threshold $t$ and reason as in the proof
of Theorem~\ref{thm:uls-sup}. Otherwise, consider the contribution of the
largest open cluster at a threshold $t_m$ such that $p_c -p_0(t_m)
\asymp m^{-\lambda/\alpha}$. As in Theorem~\ref{thm:cc-cube-cri},
the largest
open cluster will be comparable in size to, and occupy a substantial
portion of $K$. Reasoning again as in the proof of Theorem~\ref{thm:uls-sup},
the contribution is of order $m^{\alpha d/2} \theta_m \geq m^{\alpha
/\nu
} \theta_m \geq m^{\alpha/\nu-\lambda}$, which increases as a positive
power of $m$.
\section*{Appendix B: The scan statistic as the GLR}\label{sec:proof-scan}
\setcounter{equation}{0}
\renewcommand{\theequation}{B.\arabic{equation}}
We show that the simple scan statistic defined in (\ref{scan})
approximates the scan statistic of Kulldorff~\cite{Kul}, which
is strictly speaking the GLR, defined as follows. The log-likelihood
under $\bbH^m_{1, K}$ is given by
\[
\operatorname{loglik}(K, \theta, \theta_0) := |K| \bigl( \theta\bar
{X}_K -
\log
\varphi(\theta)\bigr ) + |K^c| \bigl(\theta_0 \bar{X}_{K^c} -\log
\varphi(\theta_0) \bigr).
\]
Assuming $\theta$ and $\theta_0$ are both unknown, the log GLR is
defined as
\[
\max_{K \in\Km} \sup_{\theta> \theta_0} \operatorname{loglik}(K,
\theta
, \theta
_0) -\sup_{\theta_0}\operatorname{loglik}(\Vm, \theta_0, \theta_0),
\]
which is equal to
%
\begin{equation}\label{Scan}
\max_{K \in\Km} [|K| \Lambda^*(\bar{X}_K) + |K^c| \Lambda
^*(\bar{X}_{K^c}) -|\Vm| \Lambda^*(\bar{X}_{\Vm}) ]_+.
\end{equation}
(The subscript $_+$ denotes the positive part.)

Under the normal location model, $\Lambda^*(x) = x^2/2$ and (\ref
{Scan}) is equal to
\[
\max_{K \in\Km} \frac{|\Vm| |K|}{|\Vm| - |K|} (\bar{X}_K -\bar
{X}_{\Vm})_+^2.
\]
(We used the fact that $\bar{X}_K \geq\bar{X}_{K^c} \iff\bar{X}_K
\geq\bar{X}_{\Vm}$.)
If $k_m^+ := \max\{|K|\dvt K \in\Km\}$ satisfies $k_m^+/|\Vm| \to0$,
which is the case in our examples, the fraction above is equal to $|K|
(1 + \O(k^+_m/|\Vm|))$. Moreover, knowing that there is always a
cluster $K$ such that $\bar{X}_K \geq\bar{X}_{\Vm}$, we get that the
square root of (\ref{Scan}) is approximately equal to
%
\begin{equation}\label{scan0}
\max_{K \in\Km} \sqrt{|K|} (\bar{X}_K -\bar{X}_{\Vm}),
\end{equation}
which is the version of (\ref{scan}) when $\mu_0$ is unknown. (Note
that $\bar{X}_{\Vm} = \mu_0 + \O(|\Vm|)^{-1/2}$, by the central limit
theorem, so that (\ref{scan0}) is within $\O(k^+_m/|\Vm|)^{1/2}$ from
(\ref{scan}).)
This approximation is actually valid more generally, at least in a way
that suffices for the asymptotic analysis that we perform in this work.
Indeed, with $\sigma_0^2 = \Var_0(X_v)$, we have $\Lambda^*(x) = (x
-\mu
_0)^2/(2\sigma_0^2) + \O(x -\mu_0)^3$ in the neighborhood of $\mu_0$.
Assuming that $k_m^- := \min\{|K|\dvt K \in\Km\}$ satisfies $k_m^- \to
\infty$, which is the case in our examples, the approximation of the
square root of (\ref{Scan}) by (\ref{scan0}) is valid under the null,
because $\bar{X}_K = \mu_0 + \O(k_m^-)^{-1/2}$ and $\bar{X}_{K^c},
\bar
{X}_{\Vm} = \mu_0 + \O(|\Vm|)^{-1/2}$, by the central limit theorem and
the fact that $k_m^- \to\infty$ and $k_m^+/|\Vm| \to0$. The same
applies under the alternative if $\theta_m \to0$, so that $\mu
_{\theta
_m} := \E_{\theta_m}(X_v) \to\mu_0$, and therefore, $\bar{X}_K$ for
any $K \in\Km$. When $\theta_m$ is bounded away from 0, the two
statistics, square root of (\ref{Scan}) and (\ref{scan0}), are both of
order $\sqrt{|K|}$, where $K$ denotes the cluster under the alternative
(or in the case of the $\operatorname{ULS}$  scan, the largest open cluster within the
anomalous cluster). Taken together, these findings are sufficient to
allow us to conclude that the tests based on (\ref{Scan}) and (\ref
{scan}) behave similarly.
\end{appendix}

\section*{Acknowledgements}
The authors thank Ganapati Patil for providing additional references on
the upper level set scan statistic, and Mikhail Langovoy for alerting
the authors of his work on the LOC test. They are grateful to Thierry
Bodineau and Rapha\"el Cerf for their remarks on the second-largest
cluster in supercritical percolation,
and to anonymous referees for helping them improve the presentation of
the paper.
EAC was partially supported by grants from the National Science
Foundation (DMS-06-03890) and the Office of Naval Research
(N00014-09-1-0258), as well as a Hellman Fellowship. GRG was partially
supported by the Engineering
and Physical Sciences Research Council under Grant EP/103372X/1.


\printhistory

\end{document}